\documentclass[12pt]{amsart}
\usepackage{amscd,amssymb,verbatim}
\usepackage{amssymb,amsmath,amsthm,epsfig}

\def\cnum#1{\bigcirc\kern -8pt#1}
\newcommand{\n}{\noindent}
\newcommand{\vp}{\varepsilon}

\theoremstyle{plain}
\newtheorem{thm}{Theorem}[section]

\newtheorem{cor}[thm]{Corollary}
\newtheorem{prop}[thm]{Proposition}
\newtheorem{lem}[thm]{Lemma}
\newtheorem{defn}[thm]{Definition}
\newtheorem{rem}[thm]{Remark}
\newtheorem{notation}[thm]{Notation}

\newcommand{\T}{\mathcal T}
\newcommand{\B}{\mathcal B}

\newcommand{\R}{\mathcal R}
\newcommand{\NC}{\mathcal N}
\newcommand{\supp}{{\rm supp}\,}

\newcommand{\N}{\mathbb{N}}

\newcommand{\Q}{\mathbb{Q}}

\newcommand{\mnorm}[1]{\left\vert\kern-0.9pt\left\vert\kern-0.9pt\left\vert #1
    \right\vert\kern-0.9pt\right\vert\kern-0.9pt\right\vert}

\begin{document}

\title{A Hereditarily Indecomposable Asymptotic $\ell_2$ Banach Space}

\author{G. Androulakis,  K. Beanland$^*$}
\date{}
\thanks{$^*$The present paper is part of the Ph.D thesis of the second
author which is prepared at the University of South Carolina under
the supervision of the first author}

\begin{abstract}
A Hereditarily Indecomposable asymptotic $\ell_2$ Banach space is
constructed. The existence of such a space answers a question of B.
Maurey and verifies a conjecture of W.T. Gowers.
\end{abstract}

\maketitle

\markboth{G. Androulakis,  K. Beanland}{An asymptotic $\ell_2$
Hereditarily Indecomposable Banach space}

\section{Introduction} \label{sec1}

A famous open problem in functional analysis is whether there exists
a Banach space $X$ such that every (bounded linear) operator on $X$
has the form $\lambda +K$ where $\lambda$ is a scalar and $K$
denotes a compact operator.  This problem is usually called the
``scalar-plus-compact'' problem \cite{G1}. One of the reasons this
problem has become so attractive is that by a result of N. Aronszajn
and K.T. Smith \cite{AS}, if a Banach space $X$ is a solution to the
scalar-plus-compact problem then every operator on $X$ has a
non-trivial invariant subspace and hence $X$ provides a solution to
the famous invariant subspace problem.  An important advancement in
the construction of spaces with ``few'' operators was made by W.T.
Gowers and B. Maurey \cite{GM1},\cite{GM2}. The ground breaking work
\cite{GM1} provides a construction of a space without any
unconditional basic sequence thus solving, in the negative, the long
standing unconditional basic sequence problem. The Banach space
constructed in \cite{GM1} is Hereditarily Indecomposable (HI), which
means that no (closed) infinite dimensional subspace can be
decomposed into a direct sum of two further infinite dimensional
subspaces.  It is proved in \cite{GM1} that if $X$ is a complex HI
space then every operator on $X$ can be written as $\lambda +S$
where $\lambda$ is a scalar and $S$ is strictly singular (i.e. the
restriction of $S$ on any infinite dimensional subspace of $X$ is
not an isomorphism). It is also shown in \cite{GM1} that the same
property remains true for the real HI space constructed in
\cite{GM1}.  V. Ferenczi \cite{F1} proved that if $X$ is a complex
HI space and $Y$ is an infinite dimensional subspace of $X$ then
every operator from $Y$ to $X$ can be written as $\lambda i_Y +S$
where $i_Y : Y \rightarrow X$ is the inclusion map and $S$ is
strictly singular. It was proved in \cite{GM2} that, roughly
speaking, given an algebra of operators satisfying certain
conditions, there exists a Banach space $X$ such that for every
infinite dimensional subspace $Y$, every operator from $Y$ to $X$
can be written as a strictly singular perturbation of a restriction
to $Y$ of some element of the algebra.

The construction of the first HI space prompted researchers to
construct HI spaces having additional nice properties. In other
words people tried to ``marry'' the exotic structure of the HI
spaces to the nice structure of classical Banach spaces.  The
reasons behind these efforts were twofold! Firstly, by producing
more examples of HI spaces having additional well understood
properties we can better understand how the HI property effects
other behaviors of the space. Secondly, there is hope that endowing
an HI space with additional nice properties could cause the
strictly singular and compact operators on the space to coincide
giving a solution to the scalar-plus-compact problem.

An open problem that has resisted the attempts of many experts is
whether there exists a weak Hilbert HI Banach space. Recall that an
infinite dimensional Banach space $X$ is a weak Hilbert Banach space
\cite{P1},\cite{P2} if there exist positive numbers $\delta, C$ such
that every finite dimensional space $E \subset X$ contains a
subspace $F\subset E$ such that $\dim F \geq \delta \dim E$, the
Banach-Mazur distance between $F$ and $\ell_2^{\dim F}$ is at most
equal to $C$ and there is a projection $P:X \rightarrow F$ with
$\|P\|\leq C$, ($\ell_2^n$ denotes the Hilbert space of dimension
$n$). Operator theory on weak Hilbert spaces has been studied in
\cite{P1},\cite{P2}.  In particular, the Fredhold alternative has
been established for weak Hilbert spaces.

Recall some standard notation: Given a Schauder basis $(e_n)$ of a Banach
space, a sequence $(x_n)$ of non-zero vectors of $Span(e_m)_m$ is called a
block basis of the $(e_i)$ if there exist successive subsets $F_1 <
F_2 < \cdots $ of $\N$, (where for $E,F \subset \N$, $E<F$ means
$\max E< \min F$), and a scalar sequence $(a_n)$ so that
$x_n=\sum_{i\in F_n} a_i e_i$ for every $n \in \N$.  We write $x_1 <
x_2 < \cdots$ whenever $(x_n)$ is a block basis of $(e_i)$.  If
$x=\sum_n a_n e_n \in Span(e_m)_m$ then define the support of $x$ by
$\supp x = \{i :a_i \not=0 \}$, and the range of $x$, $r(x)$, as the
smallest interval of integers containing $\supp x$.

Some of the efforts that have been made in order to construct HI
space possessing additional nice properties are the following.
Gowers \cite{G2} constructed an HI space which has an asymptotically
unconditional basis.  A Schauder basis $(e_n)$ is called
asymptotically unconditional if there exist a constant $C$ such that
for any positive integer $m$, and blocks $(x_i)_{i=1}^m$ of $(e_n)$
with $m \leq x_1$ (i.e $m \leq \min \supp x_1$) and for any signs
$(\vp_i)_{i=1}^m \subset \{\pm 1\}$ we have

\[ \biggl\| \sum_{i=1}^m \vp_i x_i \biggr\| \leq C
\biggl\| \sum_{i=1}^m x_i \biggr\|. \]

Maurey \cite[page 141-142]{M} asked whether there exists and asymptotic $\ell_p$ space for $1 < p
< \infty$ and Gowers conjectured existenc the of such spaces in
\cite[page 112]{G2}. Recall that a Banach space $X$ having a
Schauder basis $(e_n)$ is called asymptotic $\ell_2$ if there exists
a constant $C$ such that for every $m \in \N$ and all blocks
$(x_i)_{i=1}^m$ of $(e_n)_n$ with $m\leq x_1$ we have

\[ \frac{1}{C} \biggl(\sum_{i=1}^m \|x_i \|^2 \biggr)^\frac{1}{2} \leq \biggl\| \sum_{i=1}^m x_i \biggr\| \leq C
 \biggl(\sum_{i=1}^m \|x_i \|^2 \biggr)^\frac{1}{2}. \]

In the present paper we construct an HI Banach space which is
asymptotic $\ell_2$.  Our approach closely uses the methods and
techniques of the paper \cite{G} of I. Gasparis. The norm of our
space $X$ satisfies an upper $\ell_2$-estimate for blocks (i.e.
there exists a constant $C$ such that for all blocks $(x_i)_{i=1}^m$
of $(e_n)_n$ we have $\| \sum_{i=1}^m x_i \| \leq C (\sum_{i=1}^m
\|x_i \|^2 )^\frac{1}{2} $). In particular our result strengthens a
result of N. Dew \cite{D} who constructed an HI space which
satisfies an upper $\ell_2$-estimate (but not a lower estimate for
blocks $(x_i)_{i=1}^m$ with $m \leq x_1$).

S.A. Argyros and I. Deliyanni \cite{AD} constructed an HI space
which is asymptotic $\ell_1$. Recall that a Banach space $X$ having
a Schauder basis $(e_n)$ is called asymptotic $\ell_1$ if there
exists a positive constant $C$ such that for every $m \in \N$ and
all blocks $(x_i)_{i=1}^m$ of $(e_n)_n$ with $m\leq x_1$ we have

\[ \biggl\| \sum_{i=1}^m x_i \biggr\| \geq C
 \sum_{i=1}^m \|x_i \|. \]

Ferenczi \cite{F2} constructed a uniformly convex HI space. Argyros
and V. Felouzis \cite{AF} showed that for every $p >1 $ there exists
an HI space $X_p$ such that $\ell_p$ (or $c_0$ when $p=\infty$) is
isomorphic to a quotient of $X_p$.  In particular the dual space
$X_p^*$ is not an HI space since it contains an isomorph of $\ell_q$
(for $1/p + 1/q=1$).  Argyros and A. Tolias \cite{AT} have
constructed an HI space whose dual space is saturated with
unconditional sequences.

Finally we mention that a Banach space $Y$ such that every operator
on $Y$ can be written as $\lambda +S$ with $S$ strictly singular,
has to be indecomposable (i.e. the whole space cannot be
decomposable into the direct sum of two infinite dimensional
subspaces) but not HI. Indeed, Argyros and A. Manoussakis
\cite{AM1},\cite{AM2} have constructed such spaces $Y$ not
containing an HI subspace.

\section{The Construction of the Space $X$} \label{sec2}

In this section we construct a Banach space $X$.  We will prove in
section 3 that $X$ is asymptotic $\ell_2$ and in section 4 that $X$
is HI.  The construction makes use of the Schreier families
$S_{\xi}$ (for $\xi < \omega$) which are defined in the following
way, \cite{AA}. Set $S_0=\{\{n\}\}:n \in \mathbb{N}\} \cup
\{\emptyset \}$. After defining $S_{\xi}$ for $\xi < \omega $, set

$$S_{\xi +1}=\left\{ \cup_{i=1}^n F_i : n \in \N, n \leq
F_1 < \cdots < F_n, F_i \in S_{\xi}\right\},$$

\n (here we assume that the empty set satisfies $\emptyset < F$ and $F
< \emptyset$ for any set $F$). Important properties of the Schreier
families is that they are hereditary (i.e if $F\in S_{\xi}$ and
$G\subset F$ then $G \in S_{\xi}$), spreading (i.e  if
$(p_i)_{i=1}^n \in S_{\xi}$ and $p_i \leq q_i$ for all $i \leq n$
then $(q_i)_{i=1}^n \in S_{\xi}$), and they have the convolution
property (i.e. if $F_1 < \cdots < F_n$ are each members of
$S_{\alpha}$ such that $\{\min F_i : i \leq n \}$ belongs to
$S_{\beta}$ then $\cup_{i=1}^n F_i$ belongs to $S_{\alpha +
\beta}$).  For $E_i \subset \N$ we say $(E_i)_{i=1}^n$ is
$S_{\xi}$-admissible if $E_1 < E_2 < \cdots <E_n$ and $(\min E_i)_i
\in S_{\xi}$.

Let $[\N]$ denote the collection of infinite sequences of positive
integers and for $M\in [\N]$ let $[M]$ denote the collection of
infinite sequences of elements of $M$.  Let $c_{00}$ denote the
vector space of the finitely supported scalar sequences and $(e_n)$
denote the unit vector basis of $c_{00}$.

Using Schreier families we define repeated hierarchy averages
similarly as in \cite{G}.  For $\xi < \omega$ and $M \in
[\mathbb{N}]$, we define a sequence $([\xi ]_n^M)_{n=1}^{\infty}$,
of elements of $c_{00}$ whose supports are successive subsets of
$M$, as follows:

For $\xi = 0$, let $[\xi]_{n}^M = e_{m_n}$ for all $n\in
\mathbb{N}$, where $M=(m_n)$. Assume that
$([\xi]_n^M)_{n=1}^{\infty}$ has been defined for all $M\in
[\mathbb{N}]$.  Set

$$[\xi +1]_1^M=\frac{1}{\min M}\sum_{i=1}^{\min M}[\xi]_{i}^{M}.$$

\noindent Suppose that $[\xi +1]_1^M < \cdots < [\xi + 1]_n^M$ have
been defined.  Let

$$M_n=\{m\in M : m > \max \supp[\xi +1]_n^M \} ~\mbox{and} ~ k_n= \min M_n .$$

\noindent Set

$$[\xi +1]_{n+1}^M=\frac{1}{k_n}\sum_{i=1}^{k_n} [\xi]_{i}^{M_n}. $$

\n For $x \in c_{00}$ let $(x(k))_{k\in \N}$ denote the coordinates
of $x$ with respect to $(e_k)$ (i.e. $x=\sum_k x(k)e_k$).  For $M
\in [\N]$, $\xi < \omega$ and $n\in \N$ define $(\xi)_n^M \in
c_{00}$ by $(\xi)_n^M(k)=\sqrt{[\xi]_n^M(k)}$ for all $k \in \N$.
It is proved in \cite{GL} that for every $M \in [\N]$ and $\xi <
\omega$, $\sup\{\sum_{k\in F} [\xi]_1^M(k) :F\in S_{\xi-1}\} < \xi /
\min M.$  From this it follows that

\begin{equation}
\begin{split}
&\mbox{for every $\xi < \omega$ and $\vp >0$ there exists $n \in \N$ such
that for all $M \in [\N]$ with} \\
    & \mbox{$n \leq \min M$ we have that $\sup
\biggl\{ \biggl(\sum_{k\in F} ((\xi)_1^M(k))^2 \biggr)^\frac{1}{2}
:F \in S_{\xi - 1} \biggr\} < \vp.$}
\label{eqn:AA}
\end{split}
\end{equation}

\begin{defn}
Let $(u_n)_n$ be a normalized block basis of $(e_n)_n$, $\vp >0$ and
$1\leq \xi < w$. Set $p_n = \min \supp u_n$ for all $n \in
\mathbb{N}$ and $P=(p_n)$.
\begin{enumerate}
\item An $(\vp, \xi )$ squared average of $(u_n)_n$ is any vector that
can be written in the form $\sum_{n=1}^{\infty} (\xi)_1^R(p_n) u_n$,
where $R \in [P]$ and $\sup \{ (\sum_{k\in F} ((\xi)_1^R(k))^2 )^\frac{1}{2} : F \in S_{\xi - 1} \}
< \vp$.
\item A normalized  $(\vp, \xi )$ squared average of $(u_n)_n$ is any vector $u$ of the form
$u=v/\|v\|$ where $v$ is a  $(\vp, \xi )$ squared average of
$(u_n)_n$.  In the case where $\|v\| \geq 1/2$, $u$ is called a smoothly
normalized $(\vp, \xi )$ squared average of $(u_n)_n$.
\end{enumerate}
\end{defn}

In order to define the asymptotic $\ell_2$ HI space $X$ we fix four
sequences $M=(m_i)$, $L=(\ell_i)$, $F=(f_i)$ and $N=(n_i)$ of
positive integers which are defined as follows:
 Let $M=(m_i)_{i\in \N} \in [\mathbb{N}]$ be such that $m_1 >
246$ and $m_i^2< m_{i+1}$ for all $i \in \mathbb{N}$.  Choose
$L=(l_i)_{i\in \N} \in [\mathbb{N}]$ such that and $2^{l_i}>m_i$ for
all $i \in \mathbb{N}$.  Now choose and infinite sequences
$N=(n_i)_{i\in \N \cup \{0\}}$ and $F=(f_i)_{i\in \N}$ such that
$n_0=0$, $l_j(f_j + 1) <n_j$ for all $j \in \mathbb{N}$, $f_1=1$ and
for $j \geq 2$,

\begin{equation}
f_j = \max \left\{ \sum_{1 \leq i<j}\rho_i n_i : \rho_i \in
\mathbb{N}\cup \{0\}, \prod_{1 \leq i<j}m_i^{\rho_i} < m_j^3
\right\}. \label{eqn:FF}
\end{equation}

We now define appropriate trees.

\begin{defn}
A set $\T$ is called an appropriate tree if the
following four conditions hold:
\begin{enumerate}
\item $\T$ is a finite set and each element of $\T$ (which is called a node
of $\T$) is of the form $(t_1, \ldots , t_{3n})$ where $n\in \N$,
$t_{3i-2} \in M$ for $1\leq i <n$ (these nodes are called the
$M$-entries of $(t_1, \ldots , t_{3n})$), $t_{3n-2}=0$, $t_{3i-1}$
is a finite subset of $\N$ for $1\leq i\leq n$, and $t_{3i}$ is a
rational number of absolute value at most equal to 1 for $1 \leq i
\leq n$.
\item $\T$ is partially ordered with respect to the initial segment
inclusion $\prec$, i.e. if $(t_1, \ldots , t_{3n})$, $(s_1, \ldots ,
s_{3m}) \in \T$ then $(t_1, \ldots , t_{3n}) \prec (s_1, \ldots ,
s_{3n})$ if $n<m$ and $t_i=s_i$ for $i=1, \ldots , 3n$.  For
$\alpha, \beta \in \T$ we also write $\alpha \preceq \beta$ to
denote $\alpha \prec \beta$ or $\alpha =\beta$.  For $\alpha \in \T$
the elements $\beta \in \T$ satisfying $\beta \prec \alpha$
(respectively $\alpha \prec \beta$) are called the predecessors
(resp. successors) of $\alpha$.  If $(t_1, \ldots , t_{3n}) \in \T$
then the length of $(t_1, \ldots , t_{3n})$ is denoted by $|(t_1,
\ldots , t_{3n})|$ and it is equal to $3n$.  There exists a unique
element of $\T$ which has length 3 and it is called the root of
$\T$, and it is the minimum element of $\T$ with respect to $\prec$.
Every element $\alpha \in \T$ except the root of $\T$ has a unique
immediate predecessor which is denoted by $\alpha^-$.  If $\alpha$
is the root of $\T$ set $\alpha^- = \emptyset$.  If $(t_1, \ldots ,
t_{3n})\in \T$ then $(t_1, \ldots , t_{3\ell})\in \T$ for all $1
\leq \ell \leq n$. The nodes of $\T$ without successors are called
terminal.  If $\alpha \in \T$ is non-terminal, then the set of nodes
$\beta \in \T$ with $\alpha \prec \beta$ and $|\beta|=|\alpha|+3$
are called immediate successors of $\alpha$. Also $D_{\alpha}$
denotes the set of immediate successors of $\alpha$.
\item If $\alpha \in \T$ then the last three entrees of $\alpha$ will be denoted by
$m_{\alpha}$, $I_{\alpha}$ and $\gamma_{\alpha}$ respectively.  If
$\alpha$ is the root of $\T$ then $m_{\alpha}$ (resp. $I_{\alpha}$)
is called the weight (resp. the support of $\T$ denoted by
$\supp(\T)$). If $\alpha$ is a terminal node of $\T$ then
$I_{\alpha}=\{p_\alpha\}$ for some $p_\alpha \in \N$.  If $\alpha$
is a non-terminal node of $\T$ then $I_\alpha = \cup\{I_{\beta} :
\beta \in D_\alpha \}$ and for $\beta, \delta \in D_\alpha$ with
$\beta \not= \delta$ we have either $I_\beta < I_\delta$ or
$I_\delta < I_\beta$.
\item If $\alpha \in \T$ is non-terminal and $m_{\alpha}=m_{2j}$
for some $j$, then $(I_{\beta})_{\beta \in D_{\alpha}}$ is
$S_{n_{2j}}$-admissible and $\sum_{\beta \in D_{\alpha}}
\gamma_{\beta}^2 \leq 1$.
\end{enumerate}
\label{defn:tree}
\end{defn}

\n Now set

$$G=\{\T : \T ~\mbox{is an appropriate tree}\}.$$

\noindent We make the convention that the empty tree belongs to $G$.

If $\T_1 ,\T_2 \in G$ then we write $\T_1 < \T_2$ if $\supp(\T_1)<
\supp(\T_2)$.  If $\T \in G$ and $I$ is an interval of integers then
we define the restriction of $\T$ on $I$, $\T|_{I}$, to denote the
tree resulting from $\T$ by keeping only those $\alpha \in \T$ for
which $I_\alpha \cap I \not=0$ and replacing $I_\alpha$ by $I_\alpha
\cap I$.  It is easy to see that $\T|_{I} \in G$.  It $\alpha \in
\T$ set $\T_\alpha = \{\beta \setminus \alpha^- : \beta \in \T,
\alpha\preceq \beta \}$ (for $\alpha =(t_1, \ldots ,t_{3n}) \prec
\beta = (t_1,\ldots,t_{3m})$ let $\beta \setminus \alpha =
(t_{3n+1}, \ldots , t_{3m})$). Clearly $\T_\alpha \in G$. For $\T
\in G$ and $\alpha_0$ the root of $\T$, define $-\T$ by changing
$\gamma_{\alpha_0}$ to $-\gamma_{\alpha_0}$ and keeping everything
else in $\T$ unchanged.

Define an injection

$$\sigma :\{ (\T_1 < \cdots < \T_n): n\in \mathbb{N}, \T_i \in
G ~ (i \leq n) \} \rightarrow \{ m_{2j}: j\in \mathbb{N}\}$$

\noindent such that $\sigma(\T_1 , \cdots , \T_n) > w(\T_i)$ for all
$1 < i \leq n$.

\begin{defn}
\begin{enumerate}
\item For $j \in \N$, a collection $(\T_\ell)_{\ell=1}^n \subset G$
 is called $S_j$ admissible if $(\supp \T_{\ell})_{\ell=1}^n$ is $S_j$-admissible.
\item A collection of $S_{n_{2j+1}}$-admissible trees
$(\T_\ell)_{\ell=1}^n \subset G$ is called $S_{n_{2j+1}}$-dependent
if $w(\T_1)=m_{2j_1}$ for some $j_1 \geq j+1$ and $\sigma(\T_1,
\cdots , \T_{i-1})=w(\T_i)$ for all $2 \leq i \leq n$.
\item Let $G_0 \subset G$. A collection of $S_{n_{2j+1}}$ admissible trees
$(\T_\ell)_\ell^n \subset G$ is said to admit an
$S_{n_{2j+1}}$-dependent extension in $G_0$ if there exist $k \in
\N\cup \{0\}$ , $L \in \N$ and $\R_1 < \cdots < \R_{k+1} < \cdots <
\R_{k+n} \in G_0$, where $\R_{k+i}|_{[L,\infty)} = \T_i$ for all $1
\leq i\leq n$.
\item We say that $G_0 \subset G$ is self dependent if for all
$j\in \N$, $\T \in G_0$ and $\alpha \in \T$ such that
$m_{\alpha}=m_{2j+1}$, the family $\{\T_\beta: \beta \in D_\alpha\}$
admits an $S_{n_{2j+1}}$-dependent extension in $G_0$.
\end{enumerate}
\end{defn}

A set $G_0 \subset G$ is symmetric if $-\T \in G_0$ whenever $\T \in
G_0$; $G_0$ is closed under restriction to intervals if $\T|_J \in
G_0$ whenever $\T \in G_0$ and $J \subset \N$ an interval.

\begin{defn}
Let $\Gamma$ be the union of all non-empty, self-dependent,
symmetric subsets of $G$ closed under restrictions to intervals such
that for every $\T \in \Gamma$ and $\alpha \in \T$, if
$m_\alpha=m_{2j+1}$ for some $j\in \N$ then:
\begin{enumerate}
\item $\alpha$ is non-terminal.
\item $D_\alpha = \{\beta_1,\cdots , \beta_n\}$ for some $n\in \N$ with
$\T_{\beta_1} <\cdots < \T_{\beta_n}$ and there exist $k \in \N\cup
\{0\}$, $L \in \N$ and an $S_{n_{2j+1}}$-dependent family $\R_1 <
\cdots < \R_{k+1} < \cdots < \R_{k+n} \in \Gamma$, with
$\R_{k+i}|_{[L,\infty)} = \T_i$ for $1 \leq i\leq n$.
\item $(\gamma_i)_{i=1}^n$ is a non-increasing sequence of positive rationals such that
$\sum_{i=1}^n\gamma_i^2 \leq 1$.
\end{enumerate}
\label{defn:dep}
\end{defn}

\begin{notation}
Let $\T \in \Gamma$ and $\alpha \in \T$.
\begin{enumerate}
\item Define the height of the tree $\T$ by $o(\T)=\max\{|\beta|:\beta \in \T\}$.
\item Let $m(\alpha)=\Pi_{\beta \prec \alpha} m_{\beta}$ if $|\alpha| >3$, while $m(\alpha)=1$ if $|\alpha|=3$.
\item If $m_\alpha =m_i$ for some $i \in\N$ set $n_\alpha = n_i$.  Also set
$n(\alpha)=\sum_{\beta \prec \alpha} n_{\beta}$ if $|\alpha| >3$, while $n(\alpha)=0$ if $|\alpha|=3$.
\item Let $\gamma(\alpha)=\Pi_{\beta \prec \alpha} \gamma_{\beta}$ if $|\alpha| >3$, while
$\gamma(\alpha)=\gamma_{\alpha}$ if  $|\alpha|=3$.
\end{enumerate}
\end{notation}

Let $(e_n^*)_n$ denote the biorthogonal functionals to the unit
vector basis of $c_{00}$. Given $\T \in \Gamma $, set

$$x^*_{\T}=\sum_{\alpha \in \max \T} \frac{\gamma(\alpha)\gamma_{\alpha}}{m(\alpha)}e^*_{p_{\alpha}}$$

\noindent where $\max \T$ is the set of terminal nodes of $\T$ and
$I_{\alpha}= \{p_{\alpha} \}$ for $\alpha \in \max \T$.

Let $\mathcal{N}=\{x^*_{\T}:\T \in \Gamma \}$  and define $X$ to be
completion of $c_{00}$ under the norm $\|x\|=\sup\{|x^*(x)|: x^* \in
\NC \}$.

Note that for each $\T \in \Gamma$ there is a unique norming
functional $x^*_{\T}\in \NC \subset \{x^*: \|x^*\|\leq 1\}$ thus set
$w(x^*_\T)=w(\T)$ and $\supp(x^*_\T)=\supp(\T)$.  We will often use
the range of $x^* \in \NC$, $r(x^*)$, which is the smallest interval
containing $\supp(x^*)$.  If $x^* \in \NC$ and $I$ is an interval of
integers define the restriction of $x^*$ on $I$, $x^*|_I$, by
$x^*|_I(e_i)=x^*(e_i)$ if $i \in I$ and $x^*|_I(e_i)=0$ if $i
\not\in I$. It is then obvious that if $\T \in \Gamma$ and $I$ is an
interval of integers then $x^*_\T |_I = x^*_{\T|_I}$.  For $j\in \N$
and $\T_1 , \cdots , \T_n \in \Gamma$ we say that
$(x^*_{\T_{\ell}})_{\ell=1}^n$ is $S_{n_{2j+1}}$-dependent (or it
admits an $S_{n_{2j+1}}$-dependent extension) if
$(\T_{\ell})_{\ell=1}^n$ is $S_{n_{2j+1}}$-dependent (or admits an
$S_{n_{2j+1}}$-dependent extension).  Also we say that a collection
$(x_i)_{i=1}^n \subset c_{00}$ is $S_j$ admissible if $(\supp
x_i)_{i=1}^n$ is $S_j$ admissible.

The maximality of $\Gamma$ implies the following:

\begin{rem}
\begin{enumerate}
\item $e^*_n \in \NC$ for all $n \in \N$.
\item For each $\T \in \Gamma$ and $\alpha \in \T$ the tree
$\T_\alpha = \{\beta \setminus \alpha^{-} : \beta \in \T, \alpha
\preceq \beta\}$ is in $\Gamma$.
\item For every $x^* \in \NC$ with $w(x^*)=m_{2j}$ for a some $j \in
\N$ we can write
$$x^*=\frac{1}{m_{2j}}\sum_\ell \gamma_\ell x^*_\ell$$
for some $(\gamma_\ell)_\ell$ in $c_{00}$ with $\sum_\ell \gamma^2_\ell \leq 1$ and
$x^*_\ell \in \NC$ where $(\supp x^*_\ell)_\ell$ is $S_{n_{2j}}$-admissible.

\item For every $x^* \in \NC$ with $w(x^*)=m_{2j+1}$ for a some $j \in
\N$ we can write
$$x^*=\frac{1}{m_{2j+1}} \sum_{\ell =1}^t \gamma_\ell x^*_{\ell}$$

\n for some positive decreasing $(\gamma_\ell)_\ell$ in $c_{00}$
with $\sum_\ell \gamma^2_\ell \leq 1$.   Furthermore there exists $k
\in \N \cup\{0\}$,$L \in \N$ and $y^*_1 < \cdots < y^*_{k+1} <
\cdots <y^*_{k+t} $ where $(y^*_\ell)_\ell$ is
$S_{n_{2j+1}}$-dependent and $y^*_{k+i}|_{[L,\infty)}=x^*_i$ for all
$i\leq t$.
\end{enumerate}
\end{rem}

\section{Preliminary Estimates}

In this section we make some estimates similar to those in \cite{G}
that will be important in the proof that $X$ is H.I. First we show
that $X$ is asymptotic $\ell_2$.  Obviously $(e_n)_n$ is a
bimonotone unit vector basis for $X$ since the linear span of
$(e_n)_n$ is dense in $X$ and for finite intervals $I,J$ of integers
with $I \subset J$ and scalars $(a_n)_n$ we have $\|\sum_{n\in I}
a_n e_n\| \leq \|\sum_{n\in J} a_n e_n \|$ (this follows from the
fact that $\Gamma$ is closed under restrictions to intervals).

We now introduce a short remark.

\begin{rem}
If $x \in Span(e_n)_n$, $x^* =\frac{1}{m_i}\sum_j \gamma_j x_j^* \in
\NC$, $J=\{j:r(x_j^*) \cap r(x)\not= \emptyset\}$ then

$$\frac{1}{m_i} \sum_j \gamma_j x^*_j (x) \leq \biggl(\sum_{j \in J}
\gamma_j^2 \biggr)^\frac{1}{2} \|x\|.$$

\label{rem:trivial}
\end{rem}

\begin{proof}
Indeed there exists an arbitrarily small $\eta >0$ such that,

$$\biggl( \sum_{m \in J} \gamma_m^2 \biggr)^\frac{1}{2} + \eta \in \Q.$$

For $j \in J$ let $\beta_j = \gamma_j/(( \sum_{m \in J} \gamma_m^2
)^\frac{1}{2} + \eta ) \in \Q$.  Notice that $(\sum_j \beta^2_j
)^\frac{1}{2} \leq 1$ and if $i$ is odd then we have that
$(\beta_j)_j$ are non-decreasing and positive (since $(\gamma_j)_j$
are) and $(x^*_j)_{j \in J}$ has a dependent extension.  Thus $1/m_i
\sum_{j \in J} \beta_j x^*_j \in \NC$, hence

\begin{equation*}
\begin{split}
\frac{1}{m_i}\sum_j \gamma_j x_j^* (x) & = \biggl( \biggl(\sum_{j
\in J}\gamma_j^2\biggr)^{\frac{1}{2}} + \eta \biggr)
\frac{1}{m_i}\sum_{j \in J} \beta_j x_j^* (x) \\
    & \leq \biggl(\biggl(\sum_{j \in
    J}\gamma_j^2\biggr)^\frac{1}{2}+\eta\biggr)
    \|x\|
\end{split}
\end{equation*}

\n Since $\eta >0$ is arbitrary, the result follows.
\end{proof}

The next proposition shows that the norm of $X$ satisfies an upper
$\ell_2$-estimate for blocks.

\begin{prop}
If $(x_i)_{i=1}^m$ is a normalized block basis of $X$ then for any
sequence of scalars $(a_i)_i$ the following holds:

$$\biggl\|\sum_{i=1}^m a_i x_i \biggr\| \leq \sqrt{3} \biggl( \sum_{i=1}^m |a_i|^2\biggr)^{\frac{1}{2}}. $$
\label{prop:ul2}
\end{prop}

\begin{proof}

For the purposes of this proposition define $\Gamma_n = \{\T \in
\Gamma : o(\T)\leq 3n\}$ and $\mathcal{N}_n=\{x^*_\T : \T \in
\Gamma_n\}$. For $x \in c_{00}$ define $\|x\|_n=\sup \{x^*(x):x^*
\in \mathcal{N}_n \}$. Notice that the norm $\| \cdot \|$ of $X$
satisfies $\lim_{n\rightarrow \infty} \|x\|_n=\|x\|$. We will use
induction on $n$ to verify that $\| \cdot \|_n$ satisfies the statement
of the proposition.

For $n=1$, $\mathcal{N}_n=\{\gamma e_m^* : \gamma \in \Q ,
|\gamma|\leq 1,~  m \in \N\}$, and so the claim is trivial.

For the inductive step, let $x^*=1/m_k \sum_j \gamma_j x_j^* \in
\mathcal{N}_{n+1}$ (where $(x_j^*)_j \subset \NC_{n}$ is $S_{n_k}$
admissible and $\sum \gamma_i^2 \leq 1$) and let

$$Q(1)=\{1\leq i \leq m:~\mbox{there is exactly one}~ j ~\mbox{such that}~ r(x_j^*)\cap r(x_i) \not= \emptyset \},$$

\n and $Q(2)=\{1,\ldots , m\}\setminus Q(1)$.  Now apply the
functional $x^*$ to $\sum_{i=1}^m a_i x_i$ to obtain:

\begin{equation}
\begin{split}
\frac{1}{m_k}\sum_{j} \gamma_j x^*_j \biggl(\sum_{i=1}^m a_i x_i
\biggr) & =\frac{1}{m_k}\sum_{j} \gamma_j x^*_j \sum_{\substack{i\in
Q(1)\\ r(x_j^*)\cap r(x_i) \not= \emptyset}} a_i x_i +
\frac{1}{m_k}\sum_{j} \gamma_j x^*_j \sum_{i \in Q(2)} a_i x_i \\
    & \leq \frac{\sqrt{3}}{m_k}\sum_j \gamma_j  \biggl(\sum_{\substack{i\in Q(1) \\ r(x_j^*)\cap r(x_i) \not= \emptyset}} a_i^2\biggr)^{\frac{1}{2}}
    + \sum_{i \in Q(2)} a_i \frac{1}{m_k}\sum_j \gamma_j x^*_j(x_i) ,
    \label{eqn:CC}
\end{split}
\end{equation}

\n by applying the induction hypothesis for $\sum_{\{i \in Q(1):
r(x^*_j) \cap r(x_i)\}} a_i x_i$.  The above estimate continues as
follows,

\begin{equation}
\begin{split}
 &\leq \sum_j \gamma_j \biggl(\sum_{\substack{i\in Q(1) \\ r(x_j^*)\cap r(x_i) \not= \emptyset}} a_i^2 \biggr)^{\frac{1}{2}} +
\sum_{i \in Q(2)} a_i \biggl(\sum_{\{j : r(x_j^*)\cap r(x_i) \not= \emptyset \}} \gamma_j^2\biggr)^{\frac{1}{2}} \\
    & \leq \biggl(\sum_j\gamma_j^2\biggr)^\frac{1}{2} \biggl(\sum_j \sum_{\substack{i\in Q(1) \\ r(x_j^*)\cap r(x_i) \not= \emptyset}} a_i^2\biggr)^{\frac{1}{2}}
    + \sum_{i \in Q(2)} a_i \biggl(2 \sum_{j} \gamma_j^2\biggr)  \leq \sqrt{3}
    \biggl(\sum_{i}^m a_i^2\biggr)^{\frac{1}{2}},
    \label{eqn:BB}
\end{split}
\end{equation}

\n where for the first inequality of (\ref{eqn:BB}) we used that
$\sqrt{3}<m_1$ and Remark \ref{rem:trivial} and for the second
inequality of (\ref{eqn:BB}) we used the Cauchy-Schwartz inequality
and the fact that for each $j$  there are at most two values of $i
\in Q(2)$ such that $r(x^*_j) \cap r(x_i) \not= \emptyset$.  For the
third inequality of (\ref{eqn:BB}) we used $\sum_\ell \gamma_\ell^2
\leq 1$.  Combine (\ref{eqn:CC}) and (\ref{eqn:BB}) to finish the
inductive step.

\end{proof}

\begin{cor}
Let $(x_i)_{i=1}^n$ be a block basis of $X$ with $n \leq x_1$.
Then for any sequence of scalars $(a_i)_i$ the following holds:

$$\frac{1}{m_2}\biggl( \sum_{i=1}^n |a_i|^2\biggr)^{\frac{1}{2}} \leq
 \biggl\|\sum_{i=1}^n a_i x_i \biggr\| \leq \sqrt{3} \biggl( \sum_{i=1}^n |a_i|^2)\biggr)^{\frac{1}{2}}. $$
\label{cor:al2}
\end{cor}

\begin{proof}
Let $(x_i)_{i=1}^n$ be a normalized block sequence of $(e_n)$ such
that $n \leq x_1 < \cdots < x_n $ and scalars $(a_i)_{i=1}^n$. The
upper inequality follows from Proposition \ref{prop:ul2}.  Note
that, $(x_i)_{i=1}^n$ is $S_1$ admissible hence $S_{n_2}$
admissible. Find norm one functionals $(x_i^*)_{i=1}^n$ such that
$r( x^*_i) \subset r(x_i)$ and $x^*_i(x_i)=1$ for all $i\leq n$. To
establish the lower inequality apply the functional

$$\frac{1}{m_2}\sum_{i=1}^n \biggl( a_i \biggl/ \displaystyle \biggl(\sum_{i=1}^n a_i^2\biggr)^\frac{1}{2} \biggr) x^*_i$$

\n (whose norm is at most equal to one) to $\sum_{i=1}^n a_i x_i$.

\end{proof}

The next lemma is a variation of the decomposition lemma found in
\cite{G} and will be used in the proof of Lemma \ref{lem:any} and Proposition
\ref{prop:conds}

\begin{lem}
(Decomposition Lemma) Let $x^* \in \NC$.  Let $j \in \mathbb{N}$
be such that $w(x^*)< m_j$.  Then there exists as
$S_{f_j}$-admissible collection $(x^*_\alpha)_{\alpha \in L}$ and a sequence of
scalars $(\lambda_\alpha)_{\alpha \in L}$ such that $L=\cup_{i=1}^3 L_i$
and:
\begin{enumerate}
\item $x^*=\sum_{\alpha \in L} \lambda_\alpha x^*_\alpha$.
\item $\displaystyle \biggl( \sum_{\alpha \in L_1} \lambda_\alpha^2 \biggr)^{\frac{1}{2}}\leq \frac{1}{m_j^2}$,
$\displaystyle \biggl( \sum_{\alpha \in L} \lambda_\alpha^2
\biggr)^{\frac{1}{2}}\leq \frac{1}{w(x^*)}$.
\item $w(x^*_\alpha) \geq m_j$ for $\alpha \in L_2$.
\item  For $\alpha \in L_3$ there exists  $|\gamma_\alpha | \leq 1$
and $p_\alpha \in \mathbb{N}$ such that $x^*_{\alpha}=\gamma_\alpha e^*_{p_\alpha}$.
\end{enumerate}
\label{lem:dl}
\end{lem}

\begin{proof}
Since $x^* \in \NC$ there exists $\T \in \Gamma$ such that
$x^*=x^*_\T$.  Define three pairwise disjoint sets $L_1,L_2,L_3$ of
nodes of $\T$ such that for every branch $\B$ of $\T$ (i.e. a
maximal subset of $\T$ which is totally ordered with respect to
$\prec$) there is a unique $\alpha \in \B$ with $\alpha \in
\cup_{i=1}^3 L_i$. For every branch $\B$ of $\T$ choose a node
$\alpha \in \B$ which is maximal with respect to $\prec$ such that
$m(\alpha)<m_j^2$ and all $M$-entries of $\alpha^-$ are less than
$m_j$.  If $\alpha$ is non-terminal and $m_\alpha < m_j$ then
$\alpha^+\in L_1$, where $\alpha^+$ is the unique $\prec$-immediate
successor of $\alpha$ in $\B$.  Thus for $\alpha^+ \in L_1$ we have
$m_j^2\leq m(\alpha^+)<m_j^3$.  If $\alpha$ is non-terminal and
$m_\alpha \geq m_j$ then $\alpha \in L_2$.  If $\alpha$ is terminal
then $\alpha \in L_3$. Let $L=\cup_{i=1}^3 L_i$. For $\alpha \in L$
let

$$x^*_\alpha = x^*_{T_\alpha} = \frac{m(\alpha)}{\gamma(\alpha)}x^*_{\T}|_{I_\alpha} ~~ \mbox{and} ~~
\lambda_\alpha=\frac{\gamma(\alpha)}{m(\alpha)}.$$

\n Since $m_j^2 \leq m(\alpha)$ for $\alpha \in L_1$ we have

$$\biggl(\sum_{\alpha \in L_1} \lambda_\alpha^2 \biggr)^\frac{1}{2} =
\biggl(\sum_{\alpha \in L_1}
\biggl(\frac{\gamma(\alpha)}{m(\alpha)}\biggr)^2\biggr)^\frac{1}{2}
\leq \frac{1}{m_j^2} \biggl(\sum_{\alpha \in L_1} \gamma(\alpha)^2
\biggr)^\frac{1}{2} \leq \frac{1}{m_j^2},$$

\n where the last inequality follows from Definition \ref{defn:tree}
(4) and Definition \ref{defn:dep} (3).  If $\alpha \in L_2$ then
$w(\T_\alpha) \geq m_j$. If $\alpha \in L_3$ then
$x^*_{\T_\alpha}=\gamma_\alpha e^*_{p_\alpha}$. Finally, since
$m(\alpha)\geq w(x^*)$ for all $\alpha \in L$ we have,

$$\biggl(\sum_{\alpha \in L} \lambda_\alpha^2 \biggr)^\frac{1}{2} =
\biggl(\sum_{\alpha \in L} \biggl(\frac{\gamma(\alpha)}{m(\alpha)}\biggr)^2\biggr)^\frac{1}{2} \leq \frac{1}{w(x^*)}
\biggl(\sum_{\alpha \in L} \gamma(\alpha)^2 \biggr)^\frac{1}{2} \leq \frac{1}{w(x^*)}.$$

\n By applying the following Remark \ref{rem:nodes} (which also
appears in \cite{G})
 to the set $\{ I_\alpha : \alpha \in L \}$ we conclude
that $(x_\alpha^*)_{\alpha \in L}$ is $S_p$-admissible where $p=\max
\{ n(\alpha) : \alpha \in L \} \leq f_j$ by (\ref{eqn:AA}). Thus
$(x^*_{\alpha})_{\alpha \in L}$ is $S_{f_j}$ admissible.

\end{proof}

\begin{rem}
Let $\T \in \Gamma$. Let $F$ be a subset of $\T$ consisting of
pairwise incomparable nodes.  Then $\{ I_{\alpha} : \alpha \in F \}$
is $S_p$-admissible, where $p=\max \{n(\alpha):\alpha \in F\}$.
\label{rem:nodes}
\end{rem}

\begin{proof}
Proceed by induction on $o(\T)$.  For $o(\T)=3$ the assertion is
trivial. Assume the claim for all $\T \in \Gamma$ such that
$o(\T)<3n$.  Let $\T$ such that $o(\T)=3n$ and $w(\T)=m_i$.  If
$|F|=1$ the assertion is trivial, thus assume $|F|>1$.  Let
$\alpha_0$ be the root of $\T$.  Thus for all $\beta \in
D_{\alpha_0}$ the claim holds for $\T_{\beta}$.  For each $\beta \in
D_{\alpha_0}$ define, $F_{\beta}=\{\alpha\setminus \alpha_0 : \alpha
\in F, \beta \preceq \alpha \} \subset \T_\beta.$

\n We know that for every $\beta \in D_{\alpha_0}$ we have that
$\{I_{\alpha} : \alpha \in F_{\beta}\}$ is $S_{p_\beta}$ admissible
where $p_{\beta}= \max\{n_\beta(\alpha): \alpha \in F_\beta\}$ and
for every $\alpha \in \T_{\beta}$,

$$n_\beta (\alpha) = \sum_{\substack{\gamma \in \T_\beta \\ \gamma \prec \alpha \setminus \alpha_0 }} n_\gamma
=\sum_{\gamma \in \T, \gamma \prec \alpha} n_\alpha -n_i =n(\alpha) - n_i.$$

\n Thus $\{I_{\alpha} : \alpha \in F_{\beta}\}$ is
$S_{n(\alpha)-n_i}$ admissible for all $\beta \in D_{\alpha_0}$.
Also $\{I_\beta : \beta \in D_{\alpha_0} \}$ is $S_{n_i}$ admissible
so we use the convolution property of Schreier families to conclude
that $\{I_\alpha : \alpha \in F\}$ is $S_p$ admissible.

\end{proof}

\begin{lem}
Let $(u_n)$ be a normalized block basis of $(e_n)$.  Let $j \in 2 \mathbb{N}$
and let $y$ be an  $(\vp,f_j+1)$
squared average of $(u_n)$ with $\vp < 1/m_j$.  Let $(x^*_\ell)_{\ell} \in \NC$ be $S_{\xi}$-admissible, $\xi \leq f_j$  and
$(\gamma_{\ell})_{\ell}$ in $c_{00}$.
Then,

$$\sum_{\ell} \gamma_{\ell} x^*_{\ell}(y) \leq 5
\biggl(\sum_{\{\ell: r(x^*_\ell) \cap r (y) \not= \emptyset \}}
\gamma_\ell^2\biggr)^\frac{1}{2} .$$ \label{lem:fe}
\end{lem}

\begin{proof}
Let $p_n= \min \supp u_n$, $R \in [(p_n)]$ and $y=\sum_n (f_j+1)_1^R(p_n)u_n$.  Define

$$Q(1)=\{i:~\mbox{there is exactly one}~ \ell ~\mbox{such that}~ r(x^*_\ell) \cap r (u_n) \not= \emptyset \} ~\mbox{and}~$$
$$Q(2)=\{i:~\mbox{there are at least two $\ell $'s such that}~ r(x^*_\ell) \cap r (u_n) \not= \emptyset  \}.$$

\begin{align}
\biggl( \sum_\ell \gamma_\ell x^*_\ell \biggr) \biggl( \sum_{n \in Q(1)} (f_j+1)_1^R(p_n)u_n \biggr) &
\leq \sum_\ell |\gamma_\ell | \biggl| x^*_\ell
\biggl( \sum_{\substack{n\in Q(1) \\ r(x^*_\ell) \cap r (u_n) \not= \emptyset  }}(f_j+1)_1^R(p_n)u_n \biggr) \biggr| \nonumber  \\
    &\leq \sum_{\{\ell: r(x^*_\ell) \cap r(y) \not= \emptyset \}} |\gamma_\ell| \sqrt{3}\biggl(\sum_{\substack{n\in Q(1) \\ r(x^*_\ell) \cap r (u_n) \not= \emptyset
    }}((f_j+1)_1^R(p_n))^2\biggr)^{\frac{1}{2}} \label{eqn:A} \\
    & \leq \sqrt{3} \biggl( \sum_{\{\ell: r(x^*_\ell) \cap r (y) \not= \emptyset  \}}
    \gamma_\ell^2 \biggr)^\frac{1}{2} \nonumber ,
\end{align}

\n where for the second inequality we used Proposition
\ref{prop:ul2} and for the third inequality we used the
Cauchy-Schwartz inequality.  For $n$'s in $Q(2)$,

\begin{align}
\biggl| \sum_\ell \gamma_\ell x^*_\ell \biggl( \sum_{n\in
Q(2)}(f_j+1)_1^R(p_n)u_n \biggr| & \leq \sum_{n \in Q(2)}
(f_j+1)_1^R(p_n)\biggr)
\biggl| \frac{1}{m_j} \sum_{\{\ell: r(x^*_\ell) \cap r (u_n) \not= \emptyset  \}} \gamma_\ell  x_\ell^* (u_n)\biggr| m_j \nonumber \\
    & \leq \sum_{n \in Q(2)} (f_j+1)_1^R(p_n)
    \biggl(\sum_{\{\ell: r(x^*_\ell) \cap r (u_n) \not= \emptyset  \}} \gamma_\ell^2 \biggr)^{\frac{1}{2}} m_j \label{eqn:B} \\
    & \leq \biggl(\sum_{n \in Q(2)} ((f_j+1)_1^R (p_n))^2 \biggr)^{\frac{1}{2}} \biggl(\sum_{n\in Q(2)}
    \sum_{\{\ell: r(x^*_\ell) \cap r (u_n) \not= \emptyset  \}} \gamma_\ell^2 \biggr)^\frac{1}{2} m_j \nonumber \\
    &\leq  2 \biggl( 2 \sum_{\{\ell: r(x^*_\ell) \cap r (y) \not= \emptyset  \}} \gamma_\ell^2
    \biggr)^\frac{1}{2}, \nonumber
\end{align}

\n  where for the second inequality we used Remark \ref{rem:trivial}
and that $j$ is even.  For the third inequality we used the fact
that $(p_n)_{n \in Q(2)} \in 2S_{\xi}$ (i.e. the union of two sets
each which belongs to $S_{\xi}$), $\xi < f_j$, $\vp < 1/m_j$ and the
fact that for every $\ell$ there are at most two values of $n \in
Q(2)$ such that $r(x^*_\ell) \cap r(u_n) \not= \emptyset$. Combining
(\ref{eqn:A}) and (\ref{eqn:B}) we obtain the desired result since
$2\sqrt{2}+\sqrt{3}<5$.

\end{proof}

\begin{lem}
Let $(u_n)$ be a normalized block basis of $(e_n)$.  Let $\vp >0$
and $j$ be an even integer.  Then there a exists smoothly normalized
$(\vp, f_j+1)$ squared average of $(u_n)$.
\label{lem:sn}
\end{lem}

\begin{proof}
Let $P=(p_n)$ such that $p_n=\min \supp u_n$ for all $n \in \N$.  By
(\ref{eqn:AA}) assume without loss of generality that for all $R \in
[P]$, $\sup\{\sum_{k \in F}((f_j+1)_1^R(k))^2)^\frac{1}{2}:F\in
S_{f_j} \} <\vp$. Suppose that the claim is false.  For $1 \leq r \leq
l_j$ construct normalized block bases $(u_i^r)_i$ of $(u_n)$ as
follows: Set

$$u_i^1=\sum_n(f_j+1)_i^P (p_n)u_n.$$

\n It must be the case that $\|u_i^1\|<1/2$ for all $i \in \N$.
Now for each $1 <r \leq \ell_j$, if $(u_n^{r-1})_n$ has been defined let
$p_i^{r-1}=\min \supp u_i^{r-1}$, $P_{r-1}=(p_i^{r-1})$ and

$$u_i^r=\sum_n(f_j+1)_i^{P_{r-1}} (p^{r-1}_n)\frac{u^{r-1}_n}{\|u^{r-1}_n\|}.$$

\n For all $r$ and $i$, $\|u_i^r\|<1/2$.  Write
$u_1^{\ell_j}=\sum_{n\in F} a_n u_n$ for some finite set $F\subset
\N$ and $a_n >0$ with $(u_n)_{n\in F}$ being $S_{(f_j+1)\ell_j}$
-admissible and $(\sum_{n \in F} a_n^2)^\frac{1}{2} \geq
2^{\ell_j-1}$.  For $n \in F$ let $u_n^* \in X^*$,
$\|u_n^*\|=u_n^*(u_n)=1$ and $\supp u_n^* \subset r(u_n)$.  Set

$$x^*=\frac{1}{m_j}\sum_{n \in F} \biggl( a_n \biggl/ \biggl(\sum_{m\in F} a_m^2 \biggr)^\frac{1}{2} \biggr) u^*_n.$$

\n Since $(f_j+1)\ell_j < n_j$ and $j$ is even, we have that $\|x^*\|\leq 1$.  Thus

$$\frac{1}{2} > \|u_1^{\ell_j}\|\geq x^*(u_1^{\ell_j}) = \frac{1}{m_j}
\sum_{n \in F} \frac{a_n}{\biggl(\sum_{m\in F} a_m^2
\biggr)^{\frac{1}{2}} } u^*_n \biggl(\sum_{n\in F} a_n u_n \biggr)
\geq \frac{2^{\ell_j-1}}{m_j},$$

\n contradicting that $m_j\leq 2^{\ell_j}$.

\end{proof}

\begin{lem}
Let $(u_n)_n$ be a normalized block basis of $(e_n)_n$ and $j_0 \in \N$.  Suppose that
$(y_{k})_k$ is a block basis of $(u_n)_n$ so that $y_{k}$ is a smoothly
normalized  $(\vp_{k},f_{2j_k}+1)$ squared average of $(u_n)_n$ with
$\vp_{k}<1/m_{2j_k}$ and $j_0 < 2j_1 < 2j_2 < \cdots $.  Let $(x_m^*)_m \subset \NC$
be $S_{\xi}$ admissible, $\xi < n_{j_0}$ and $(\gamma_m)_m , (\beta_k)_k \in c_{00}$. Then

$$\sum_m \gamma_m x^*_m \biggl( \sum_k \beta_k y_{k} \biggr) \leq 22 \biggl( \sum_m
\gamma_m^2 \biggr)^\frac{1}{2} \biggl( \sum_k \beta^2_k \biggr)^\frac{1}{2}.$$
\label{lem:el}
\end{lem}

\begin{proof}
Define the following two sets,

$$Q(1)=\{k:~\mbox{there is exactly one}~ m ~\mbox{such that}~ r(x_m^*)\cap r(y_{k}) \not= \emptyset \},$$
$$Q(2)=\{k:~\mbox{there are at least two $m$'s such that}~ r(x_m^*)\cap r(y_{k}) \not= \emptyset \}.$$

\begin{equation*}
\begin{split}
\biggl|\sum_m \gamma_m & x^*_m \biggl( \sum_{\{k \in Q(1) :
r(x_m^*)\cap r(y_{k}) \not= \emptyset\}} \beta_k y_{k} \biggr)
\biggr| +
\biggl|\sum_m  \gamma_m  x^*_m \biggl( \sum_{k \in Q(2)}  \beta_k y_{k} \biggr)\biggr| \\
    & \leq \sum_m | \gamma_m | \sqrt{3} \biggl( \sum_{\{k \in Q(1) : r(x_m^*)\cap r(y_{k}) \not= \emptyset\}}  \beta_k^2 \biggr)^\frac{1}{2}
    +  \sum_{k \in Q(2)} |\beta_k | \biggl| \sum_{\{m: r(x_m^*)\cap r(y_{k}) \not= \emptyset \}} \gamma_m x^*_m (y_{k}) \biggr| \\
    & \leq 2 \biggl( \sum_m \gamma_m^2 \biggr)^\frac{1}{2} \biggl( \sum_m \sum_{\{k\in Q(1): r(x_m^*)\cap r(y_{k}) \not= \emptyset\}}
    \beta_k^2 \biggr)^\frac{1}{2} + 10 \sum_{k \in Q(2)} |\beta_k |
    \biggl( \sum_{\{m: r(x_m^*)\cap r(y_{k}) \not= \emptyset \}} \gamma_m^2 \biggr)^\frac{1}{2} \\
    & \leq 2 \biggl( \sum_m \gamma_m^2 \biggr)^\frac{1}{2} \biggl( \sum_k \beta_k^2 \biggr)^\frac{1}{2}
    + 20 \biggl( \sum_m \gamma_m^2 \biggr)^\frac{1}{2} \biggl( \sum_k \beta^2_k \biggr)^\frac{1}{2}
    \leq 22 \biggl( \sum_k \beta^2_k \biggr)^\frac{1}{2}.
\end{split}
\end{equation*}

\n For the first inequality we used Proposition \ref{prop:ul2}.  For
the second inequality we applied the Cauchy-Schwartz inequality in
the first term of the sum and used the fact that $\xi < n_{j_0} <
f_{2j_k}$ for all $k$ to apply Lemma \ref{lem:fe} in the second term
of the sum.   The ``10'' in the second part of the second inequality
comes from the fact that $y_{k}$ is smoothly normalized.  For the
third inequality we used the Cauchy-Schwartz inequality. The ``20''
after the third inequality comes from the fact that for every $m$
there are at most two values of $k \in Q(2)$ such that $r(x_m^*)\cap
r(y_{k}) \not= \emptyset $.
\end{proof}

\begin{lem}
Let $(u_n)_n$ be a normalized block basis of $(e_n)_n$.  Suppose
that $(y_j)_j$ is a block basis of $(u_n)_n$ so that $y_j$ is a
smoothly normalized  $(\vp_j,f_{2j}+1)$ squared average of $(u_n)_n$
with $\vp_j<1/m_{2j}$.  Then there exists a subsequence $(y_j)_{j\in
I}$ of $(y_j)_j$  such that for every $j_0 \in \N$, $j_1,j_2,\ldots
\in I$ with $j_0< 2j_1 <2j_2< \ldots $, $x^* \in \NC$ with
$w(x^*)\geq m_{j_0}$ and scalars $(\beta_j)_{j} \in c_{00}$ we have
that:
\begin{enumerate}
\item If $ w(x^*) < m_{j_1}$ then

$$x^*\biggl(\sum_{k} \beta_k y_{j_k} \biggr)< \frac{5}{m_e}
\biggl(\sum_{\{k: r(x^*)\cap r(y_{j_k}) \not= \emptyset \}}
\beta_k^2 \biggr)^{\frac{1}{2}},$$

\n where $m_e= m_{j_0}$ if $w(x^*)=m_{j_0}$ and $m_e=m^2_{j_0}$ if $w(x^*)>m_{j_0}$.

\item If $m_{2j_s} \leq w(x^*) < m_{2j_{s+1}}$ for some $s \geq 1$ then

$$x^*\biggl(\sum_{k \not=s } \beta_k y_{j_k} \biggr)< \frac{5}{m^2_{j_0}}
\biggl(\sum_{\{k \not = s : r(x^*)\cap r(y_{j_k}) \not= \emptyset
\}} \beta_k^2 \biggr)^{\frac{1}{2}}.$$
\end{enumerate}

\label{lem:lwl}
\end{lem}

\begin{proof}
Lemma \ref{lem:sn} assures the existence the block sequence
$(y_j)_j$ such that each $y_j$ is a smoothly normalized
$(\vp_j,f_{2j}+1)$ squared average of $(u_n)_n$. Let $T=(t_n)$,
where  $t_n=\min \supp u_n$. Choose $I=(j'_k)_k \in [ \N ]$ such
that, $j'_1$ is an arbitrary integer, $\sum_{i>k} \vp^2_{j'_i} <
\vp_{j'_k}^2$ and

\begin{equation}
\biggl(\sum_{i<k} \|y_{j'_i}\|_{\ell_1}^2 \biggr)^\frac{1}{2}<
\frac{m_{2j'_k}}{m_{2j'_{k-1}}}. \label{eqn:EP}
\end{equation}

Let $j_0 \in \N$, $j_1,j_2,\ldots \in I$ with $j_0< 2j_1 < 2j_2<
\ldots , (\beta_k)\in c_{00}$ and $x^* \in \NC$ such that $m_{j_0}
\leq w(x^*) < m_{2j_1}$. By definition
$y_{j_k}=v_{j_k}/\|v_{j_k}\|$, where $v_{j_k}=\sum_n
(f_{2j_k}+1)_1^{R_k}(t_n) u_n$ and $R_k \in [(t_n)]$ is chosen as in
Definition 2.1.  Let $x^*= 1/m_i \sum_\ell \gamma_\ell x^*_\ell$ for
some $i$ where $\sum_\ell \gamma_\ell^2 \leq 1$, $(x_\ell^*)_\ell$
is $S_{n_i}$ admissible and $i <2j_1$. Define the following two
sets.

$$Q(1)=\{n: ~\mbox{there is exactly one} ~\ell ~\mbox{such that}~ r(x^*_\ell)\cap r(u_n) \not= \emptyset \},$$

$$Q(2)=\{n: ~\mbox{there are at least two $\ell $'s such that}~ r(x^*_\ell)\cap r(u_n) \not= \emptyset \}.$$

\n We proceed with the case $n \in Q(1)$.

\begin{equation}
\begin{split}
\frac{1}{w(x^*)} \sum_\ell \gamma_\ell x^*_\ell
& \biggl( \sum_k \frac{\beta_k}{\|v_{j_k}\|} \sum_{\{ n \in Q(1): r(x^*_\ell)\cap r(u_n) \not= \emptyset \} } (f_{2j_k}+1)_1^{R_k}(t_n) u_n \biggr) \\
    &\leq \frac{2}{w(x^*)} \sum_\ell | \gamma_\ell | \biggl| x^*_\ell \biggl( \sum_{\{ k: r(x^*_\ell)\cap r(y_{j_k}) \not= \emptyset \}}
    \beta_k \sum_{\{ n \in Q(1): r(x^*_\ell)\cap r(u_n) \not= \emptyset \} } (f_{2j_k}+1)_1^{R_k}(t_n) u_n \biggr)\biggr|\\
    & \leq \frac{2\sqrt{3}}{w(x^*)} \sum_\ell | \gamma_\ell |  \biggl( \sum_{\{ k: r(x^*_\ell)\cap r(y_{j_k}) \not= \emptyset \}}
    \sum_{\{ n \in Q(1): r(x^*_\ell)\cap r(u_n) \not= \emptyset \} } \beta_k^2 ((f_{2j_k}+1)_1^{R_k}(t_n))^2  \biggr)^\frac{1}{2} \\
    & \leq \frac{4}{w(x^*)} \biggl(\sum_{\ell }
    \gamma_\ell^2 \biggr)^{\frac{1}{2}} \biggl(\sum_{\{k : r(x^*_\ell) \cap r(v_{j_k}) \not= \emptyset\}}
    \beta_k^2 \biggr)^{\frac{1}{2}} \leq \frac{3}{w(x^*)} \biggl(
    \sum_{\{k: r(x^*)\cap r(y_{j_k}) \not= \emptyset \}} \beta_k^2
    \biggr)^\frac{1}{2}
    \label{eqn:C}
\end{split}
\end{equation}

\n where for the first and second inequalities we used Proposition
\ref{prop:ul2} and the fact that $\|v_{j_k}\| < 1/2$.  For the third
inequality we used that $2\sqrt{3} \leq  4$ and the Cauchy-Schwartz
inequality. Notice that since $w(x^*) > m_{j_0}$ then we have
$4/w(x^*)< 1/m^2_{j_0}$. For  $n \in Q(2)$ we have,

\begin{equation}
\begin{split}
\frac{1}{m_i} & \sum_\ell \gamma_\ell x^*_\ell \sum_k \beta_k \frac{\sum_{n \in Q(2)}
(f_{2j_k}+1)_1^{R_k}(t_n) u_n}{\|v_{j_k}\|} \\
    & \leq 2 \sum_{\{ k :r(x^*) \cap r(y_{j_k}) \not= \emptyset\}} |\beta_k | \sum_{n \in Q(2)} (f_{2j_k}+1)_1^{R_k}(t_n)
    \biggl| \frac{1}{m_i} \sum_\ell \gamma_\ell x^*_\ell (u_n) \biggr| \\
    & \leq 2 \sum_{\{ k :r(x^*) \cap r(y_{j_k}) \not= \emptyset\}} |\beta_k| \sum_{n \in Q(2)} (f_{2j_k}+1)_1^{R_k}(t_n)
    \biggl(\sum_{\{\ell:~ r(x^*_\ell)\cap r(u_n) \not= \emptyset \}}\gamma_\ell^2\biggr)^\frac{1}{2} \\
    & \leq 2 \sum_{\{ k :r(x^*) \cap r(y_{j_k}) \not= \emptyset\}} |\beta_k| \biggl(\sum_{n \in Q(2)}
    ((f_{2j_k}+1)_1^{R_k}(t_n))^2 \biggr)^{\frac{1}{2}}
    \biggl(\sum_{n \in Q(2)} \sum_{\{\ell:~ r(x^*_\ell)\cap r(u_n) \not= \emptyset \}} \gamma_\ell^2 \biggr)^{\frac{1}{2}},
\label{eqn:D}
\end{split}
\end{equation}

\n where for the first inequality we used that $\|v_{j_k}\| >1/2$,
for the second inequality we used Remark \ref{rem:trivial} and for
the third inequality we used the Cauchy-Schwartz inequality.  Note
that $(u_n)_{n \in Q(2)}$ is $2S_{n_i}$ admissible (i.e. it can be
written as a union of two sets each of which is $S_{n_i}$
admissible) and $n_i \leq f_{2j_1}$. Also note that for every $\ell
$ there are at most two values of $n \in Q(2)$ such that
$r(x^*_\ell)\cap r(u_n) \not= \emptyset$, to continue (\ref{eqn:D})
as follows:

\begin{equation}
\begin{split}
    & \leq 4 \sum_{\{ k : r(x^*) \cap r(y_{j_k}) \not= \emptyset \}} \beta_k \vp_{j_k}
    \biggl(\sum_{\ell } \gamma_\ell^2\biggr)^{\frac{1}{2}}
    \leq 4 \biggl(\sum_{\{k :r(x^*) \cap r(y_{j_k})\not= \emptyset\}} \beta_k^2\biggr)^\frac{1}{2}
    \biggl(\sum_{k\geq 1} \vp_{j_k}^2\biggr)^\frac{1}{2} \\
    & \leq 8 \vp_{j_1} \biggl(\sum_{\{k : r(x^*) \cap r(y_{j_k})\not= \emptyset\}} \beta_k^2\biggr)^\frac{1}{2}
    \leq \frac{4}{m_{j_1}} \biggl(\sum_{\{k : r(x^*) \cap r(y_{j_k})\not= \emptyset\}} \beta_k^2\biggr)^{\frac{1}{2}}.
    \label{eqn:E}
\end{split}
\end{equation}

\n Obviously (\ref{eqn:C}),(\ref{eqn:D}) and (\ref{eqn:E}) finish
the proof of part (1).  Assume for $s\geq 1$, $m_{2j_s} \leq w(x^*)
< m_{2j_{s+1}}$. Estimate $x^*(\sum_{k>s} \beta_k y_{j_k})$
similarly to (\ref{eqn:C}),(\ref{eqn:D}) and (\ref{eqn:E}) where we
replace $m_{2j_1}$ by $m_{2j_{s+1}}$.  Estimate $x^* (\sum_{k<s}
\beta_k y_{j_k})$ as follows:

\begin{equation*}
\begin{split}
\frac{1}{w(x^*)} \sum_\ell \gamma_\ell x^*_\ell \sum_{k<s} \beta_k y_{j_k} &
= \sum_{\{k<s : r(x^*) \cap r(y_{j_k})\not= \emptyset\}} \biggl| \beta_k \frac{1}{w(x^*)} \sum_\ell \gamma_\ell x^*_\ell (y_{j_k}) \biggr| \\
    & \leq \sum_{\{k<s : r(x^*) \cap r(y_{j_k})\not= \emptyset\}} | \beta_k | \frac{1}{w(x^*)} \|y_{j_k}\|_{\ell_1}
    ~~(\mbox{since} ~ \|\sum_\ell \gamma_\ell x^*_\ell \|_\infty \leq 1) \\
    & \leq \biggl( \sum_{\{k<s : r(x^*) \cap r(y_{j_k})\not= \emptyset\}} \beta_k^2 \biggr)^{\frac{1}{2}} \frac{1}{w(x^*)}
    \biggl(\sum_{\{k<s : r(x^*) \cap r(y_{j_k})\not= \emptyset\}} \|y_{j_k}\|_{\ell_1}^2 \biggr)^{\frac{1}{2}} \\
    & \leq \biggl( \sum_{\{k<s :r(x^*) \cap r(y_{j_k})\not= \emptyset\}} \beta_k^2 \biggr)^{\frac{1}{2}} \frac{1}{w(x^*)} \frac{m_{j_s}}{m_{j_{s-1}}} \leq
    \frac{1}{m_{j_{s-1}}} \biggl( \sum_{\{k<s : r(x^*) \cap r(y_{j_k})\not= \emptyset\}} \beta_k^2 \biggr)^{\frac{1}{2}},
\end{split}
\end{equation*}

\n where the third inequality comes from equation (\ref{eqn:EP}).  This
finishes the proof of part (2).

\end{proof}

\begin{rem}
Lemma \ref{lem:lwl} will be used several times as follows:  Given a
normalized block sequence $(u_n)$ of $(e_n)$,  Lemma \ref{lem:sn}
will guarantee the existence of a block sequence $(y_j)$ of $(e_n)$
such that $y_j$ is a smoothly normalized $(\vp_j, f_{2j}+1)$ squared
average of $(u_n)$ and $\vp_j<1/m_{2j}$ for all $j \in \N$.  Choose
a subsequence $(y_j)_{j \in I}$ of $(y_j)_{j}$ to satisfy the
conclusion of Lemma \ref{lem:lwl}.  Let $j_0 \in \N$ and $j_1, j_2,
\ldots  \in I$ with $j_0 <2j_1 <2j_2 < \ldots $.  Let $p_k = \min
\supp (y_{j_k})$ for all $k \in \N$, $R \in [(p_k)]$,

$$y= \sum_k (n_{j_0})_1^R(p_k)y_{j_k} ~~ and ~~ g=\frac{y}{\|y\|}$$

\n be a normalized $(1/m_{j_0}^2,n_{j_0})$ squared average of $(y_j)_{j \in
I}$.  Then the conclusion of Lemma \ref{lem:lwl} will be valid for
``$\beta_k$''$=(n_{j_0})_1^R(p_k)$ and for all $x^* \in \NC$ with
$w(x^*)\geq m_{j_0}$. \label{rem:barb}
\end{rem}

\begin{lem}
Let $(y_j)_{j\in I}$, $j \in 2\N$ and a normalized $(1/m^2_{j_0})$
squared average $g$ of $(y_j)_{j\in I}$ be chosen as in Remark
\ref{rem:barb}. Then for any $S_{\xi}$ admissible family
$(x^*_{\ell})_{\ell} \subset \NC$, $\xi < n_{j_0}$ where
$w(x^*_{\ell}) \geq m_{j_0}$ for all $\ell$ and
$(\gamma_{\ell})_{\ell}\in c_{00}$, we have

\begin{equation}
\begin{split}
& \sum_{\ell} \gamma_{\ell} x^*_{\ell} (g)
    < \frac{47}{m_{j_0}} \biggl(\sum_{\{\ell :~ w(x_\ell^* ) > m_{j_0}, ~r(x^*_\ell) \cap r(u) \not= \emptyset \}}
    \gamma_\ell^2 \biggr)^{\frac{1}{2}}
    + 6 \biggl(\sum_{\{\ell : ~w(x_\ell^* ) = m_{j_0}, ~r(x^*_\ell) \cap r(u) \not= \emptyset \}} \gamma_\ell^2 \biggr)^{\frac{1}{2}}.
\end{split}
\end{equation}

\label{lem:geq}
\end{lem}

\begin{proof}
Let $g=y/\|y\|$ and note that since $j_0$ is even, $\|y\| \geq
1/m_{j_0}$ where $y=\sum_k (n_{j_0})_1^R(p_k) y_{j_k}$ for $p_k =
\min \supp y_{j_k}$ for all $k \in \N$ and $R\in [(p_k)]$. Thus

\begin{equation}
\sum_\ell \gamma_\ell x^*_\ell (g) \leq m_{j_0} \sum_\ell
|\gamma_\ell || x^*_\ell (y) | \label{eqn:F}
\end{equation}

\n Let $B=\{\ell: w(x^*_\ell)> m_{j_0} \}$ and
$E=\{\ell: w(x^*_\ell)= m_{j_0} \}$.  Define,

$$Q(1)=\{k: ~\mbox{there is exactly one} ~\ell ~\mbox{such that}~ r(x^*_\ell) \cap r(y_{j_k}) \not= \emptyset \},$$

$$Q(2)=\{k: ~\mbox{there are at least two $\ell $'s such that}~ r(x^*_\ell) \cap r(y_{j_k}) \not= \emptyset\}.$$

\n For $k$'s in $Q(2)$ we have,

\begin{equation}
\begin{split}
m_{j_0} &  \sum_\ell | \gamma_\ell| \biggl| x^*_\ell \biggl( \sum_{k \in Q(2)} (n_{j_0})_1^R(p_k) y_{j_k} \biggr) \biggr|  \\
    & \leq m_{j_0} \sum_{k \in Q(2)} (n_{j_0})_1^R(p_k) \sum_{ \{ \ell : r(x^*_\ell) \cap r(y_{j_k}) \not= \emptyset \}} |\gamma_\ell |
    |x^*_\ell (y_{j_k})| \\
    & \leq m_{j_0} \sum_{k \in Q(2)} (n_{j_0})_1^R(p_k) 10 \biggl( \sum_{\{\ell : r(x^*_\ell) \cap r(y_{j_k}) \not= \emptyset \}} \gamma^2_\ell
    \biggr)^\frac{1}{2},
    \label{eqn:G}
\end{split}
\end{equation}

\n where for the second inequality we used that $y_{j_k}'s$ are
smoothly normalized,  and since, $(x^*_\ell)_\ell$ are $S_{\xi}$
admissible with $\xi < n_{j_0} < f_{2j_k}$ for all $k$, we applied
Lemma \ref{lem:fe}.  By the Cauchy-Schwartz inequality the estimate
(\ref{eqn:G}) continues as follows:

\begin{equation}
\begin{split}
    & \leq 10 m_{j_0} \biggl( \sum_{k \in Q(2)} ((n_{j_0})_1^R(p_k))^2 \biggr)^\frac{1}{2}
    \biggl( \sum_{k \in Q(2)} \sum_{\{ \ell : r(x^*_\ell) \cap r(y_{j_k}) \not= \emptyset\}}  \gamma_\ell^2 \biggr)^\frac{1}{2} \\
    & \leq 10 m_{j_0} \frac{2}{m^2_{j_0}} 2
    \biggl( \sum_{\{ \ell : r(x^*_\ell) \cap r(u) \not= \emptyset\}}  \gamma_\ell^2 \biggr)^\frac{1}{2} \\
    & \leq \frac{40}{m_{j_0}} \biggl( \sum_{\{\ell \in B : r(x^*_\ell) \cap r(u) \not= \emptyset \}}  \gamma_\ell^2 \biggr)^\frac{1}{2}
    + \biggl( \sum_{\{ \ell \in E : \supp r(x^*_\ell) \cap r(u) \not= \emptyset \}}  \gamma_\ell^2 \biggr)^\frac{1}{2},
    \label{eqn:H}
\end{split}
\end{equation}

\n  For the second inequality we used the fact that $(y_{j_k})_{k\in
Q(2)}$ is $2S_{\xi}$ admissible for $\xi < n_{j_0}$, and that for
every $\ell$ there are at most two values of $k \in Q(2)$ such that
$r(x_\ell^*) \cap r(y_{j_k})\not= \emptyset$.

For each $\ell$ let $s_\ell$ be the integer $s$ such that $m_{2j_s}
\leq w(x^*_\ell) < m_{2j_{s+1}}$ and $r(x^*_\ell) \cap r(y_{j_k})
\not= \emptyset$  if such $s$ exists (obviously, no such $s$ exists
if $\ell \in E$ i.e. is defined for certain values of $\ell \in B$).  For $k$'s in $Q(1)$,

\begin{align}
m_{j_0} \sum_{\ell} | \gamma_{\ell} | &\biggl| x^*_{\ell}
\biggl( \sum_{k \in Q(1)} (n_{j_0})_1^R(p_k) y_{j_k} \biggr) \biggr| \nonumber \\
    & \leq m_{j_0} \biggl(  \sum_\ell | \gamma_\ell |
    \biggl| x^*_\ell \biggl( \sum_{\{k \in Q(1): k \not=s_\ell\}} (n_{j_0})_1^R(p_k) y_{j_k} \biggr)
    \biggr| + \sum_\ell |\gamma_\ell | \biggl| x^*_\ell \biggl((n_{j_0})_1^R(p_{s_\ell}) y_{j_{s_\ell}} \biggr) \biggr|
    \biggr).
    \label{eqn:split}
\end{align}

\n For the first term of the sum,

\begin{align}
m_{j_0} &   \sum_\ell |\gamma_\ell| \biggl| x^*_\ell \biggl(
\sum_{\{k \in Q(1): k \not=s_\ell \}} (n_{j_0})_1^R(p_k) y_{j_k} \biggr) \biggr| \nonumber \\
    & \leq m_{j_0} \biggl(  \sum_{\ell \in B}  |\gamma_\ell |\biggl| x^*_\ell \biggl(
    \sum_{\{k \in Q(1): k \not=s_\ell \}} (n_{j_0})_1^R(p_k) y_{j_k} \biggr) \biggr|  +
     \sum_{\ell \in E} | \gamma_\ell | \biggl|x^*_\ell \biggl(
    \sum_{\{k \in Q(1): k \not=s_\ell\}} (n_{j_0})_1^R(p_k) y_{j_k} \biggr) \biggr|   \biggr) \nonumber \\
    & \leq m_{j_0} \biggl(  \sum_{\ell \in B} | \gamma_\ell | \frac{5}{m^2_{j_0}} \biggl(
    \sum_{\substack{k \in Q(1), k \not=s_\ell \\ r(x^*_\ell) \cap r(y_{j_k}) \not= \emptyset }} ((n_{j_0})_1^R(p_k))^2 \biggr)^\frac{1}{2}  +
    \sum_{\ell \in E} |\gamma_\ell | \frac{5}{m_{j_0}} \biggl(
    \sum_{\substack{k \in Q(1), k \not=s_\ell \\ r(x^*_\ell) \cap r(y_{j_k}) \not= \emptyset }}
    ((n_{j_0})_1^R(p_k))^2 \biggr)^\frac{1}{2} \biggr)  \label{eqn:II} \\
    & \leq \frac{5}{m_{j_0}} \biggl( \sum_{\ell \in B} \gamma_\ell^2 \biggr)^\frac{1}{2} +
    5 \biggl( \sum_{\ell \in E} \gamma_\ell^2 \biggr)^\frac{1}{2}. \nonumber
\end{align}

\n For the second inequality of (\ref{eqn:II}) we applied Lemma
\ref{lem:lwl}.  Notice the $s_\ell$'s where picked to coincide with
part (2) of Lemma \ref{lem:lwl}.  The final inequality
of (\ref{eqn:II}) followed from the Cauchy-Schwartz inequality.

For the second part of the right hand side of (\ref{eqn:split})
notice that the only $\ell$'s that appear are the ones for which
$s_\ell$ is defined. Also recall that if $s_\ell$ is defined then
$w(x_\ell^*) > m_{j_0}$ hence $\ell \in B$. Thus the second part of
the right hand side of (\ref{eqn:split}) can be estimated as
follows:

\begin{align}
m_{j_0}  \sum_{\{\ell : ~s_\ell~ \text{is defined}  \}} &
|\gamma_\ell| \biggl| x^*_\ell \biggl( (n_{j_0})_1^R(p_{s_\ell})
y_{j_{s_\ell}} \biggr) \biggr|
 \leq  m_{j_0} \sum_{\{\ell :~s_\ell~ \text{is defined}  \}}
|\gamma_\ell | (n_{j_0})_1^R(p_{s_\ell}) \nonumber \\
    & \leq  m_{j_0} \biggl( \sum_{\{\ell :~s_\ell~ \text{is defined}  \}} \gamma_\ell^2
    \biggr)^\frac{1}{2}  \label{eqn:JJ}
    \biggl( \sum_{\{\ell :~s_\ell~ \text{is defined}  \}} ((n_{j_0})_1^R(p_{s_\ell}))^2 \biggr)^\frac{1}{2} \\
    & \leq \frac{2}{m_{j_0}} \biggl( \sum_{\{\ell : ~s_\ell~ \text{is defined}  \}}
    \gamma_\ell^2 \biggr)^\frac{1}{2}, \nonumber
\end{align}

\n where for the second inequality we applied the Cauchy-Schwartz
inequality and for the third inequality we used that
$(x^*_\ell)_\ell$ is $S_{\xi}$ admissible for $\xi < n_{j_0}$ hence
$\{(p_{s_\ell}) : s_\ell ~\mbox{is defined}\} \in 2 S_{\xi}$. The
result follows by combining the estimates (\ref{eqn:F}),
(\ref{eqn:G}),(\ref{eqn:H}),(\ref{eqn:split}),(\ref{eqn:II}), and
(\ref{eqn:JJ}).
\end{proof}

\begin{lem}
Let $(y_j)_{j\in I}$, $j_0 \in 2\N$ and a normalized
$(1/m^2_{j_0},n_{j_0})$ squared average $g$ of $(y_j)_{j\in I}$
chosen as in Remark \ref{rem:barb}. Then for any $S_{n_i}$
admissible family $(x^*_{\ell})_{\ell} \subset \NC$, $i < j_0$ and
$(\gamma_{\ell})_{\ell}\in c_{00}$, we have

\begin{equation}
\begin{split}
& \sum_{\ell} \gamma_{\ell} x^*_{\ell} (g)
    < \frac{123}{m_{e}} \biggl(\sum_{\{\ell : w(x_\ell^* ) \not= m_{j_0}, r(x^*_\ell) \cap r(g) \not= \emptyset \}}
    \gamma_\ell^2 \biggr)^{\frac{1}{2}}
    + 6 \biggl(\sum_{\{\ell : w(x_\ell^* ) = m_{j_0},
    r(x^*_\ell) \cap r(g) \not= \emptyset \}} \gamma_\ell^2
    \biggr)^{\frac{1}{2}},
    \label{eqn:star}
\end{split}
\end{equation}

\n where $m_e=\min_\ell \{w(x_\ell^*),m_{j_0}\}$.
\label{lem:any}
\end{lem}

\begin{proof}

Let $g=y/\|y\|$ and since $j_0$ is even note that $\|y\| \geq
1/m_{j_0}$ where $y=\sum_k (n_{j_0})_1^R(p_k) y_{j_k}$ for $p_k =
\min \supp y_{j_k}$ for all $k \in \N$ and $R\in [(p_k)]$. Let
$S=\{\ell: w(x^*_\ell)< m_{j_0} \}$, $E=\{\ell: w(x^*_\ell)= m_{j_0}
\}$ and $B=\{\ell: w(x^*_\ell) >m_{j_0} \}$. Using Lemma
\ref{lem:dl} for $\ell \in S$ we write

$$x^*_\ell = \sum_{m \in L_\ell} \lambda_{\ell,m} x^*_{\ell,m}, $$

\n where $L_\ell = \cup_{i=1}^3 L_{\ell,i}$ and the following are satisfied:

\begin{equation}
 \biggl( \sum_{m \in L_{\ell,1}} \lambda_{\ell,m}^2 \biggr)^{\frac{1}{2}}\leq \frac{1}{m_{j_0}^2},
\displaystyle \biggl( \sum_{m \in L_\ell} \lambda_{\ell,m}^2
\biggr)^{\frac{1}{2}}\leq \frac{1}{w(x^*_\ell)}, \label{eqn:K}
\end{equation}

\n $w(x^*_{\ell,m}) \geq m_{j_0}~\mbox{for}~m \in L_{\ell,2},~
x^*_{\ell,m}=\gamma_{\ell,m} e^*_{p_{\ell,m}}~\mbox{for}~ m \in
L_{\ell,3}, ~|\gamma_{\ell,m} | \leq 1  ~\mbox{and} ~ p_{\ell,m} \in
\mathbb{N}$. Now we have

\begin{align}
 \sum_{\ell} \gamma_{\ell} x^*_{\ell} (g) & \leq \biggl| \sum_{\{\ell \in S :r(x^*_\ell) \cap r(g) \not= \emptyset \}} \gamma_\ell
 \sum_{m \in L_\ell} \lambda_{\ell,m} x^*_{\ell,m} (g) \biggr| + \biggl| \sum_{\ell \not\in S}
 \gamma_\ell x^*_\ell (g) \biggr| \nonumber \\
    & \leq \biggl| \sum_{\{\ell \in S : r(x^*_\ell) \cap r(g) \not= \emptyset \}} \gamma_\ell
    \sum_{m \not\in L_{\ell,2}} \lambda_{\ell,m} x^*_{\ell,m} (g) \biggr|+
    \biggl| \sum_{\{\ell \in S :r(x^*_\ell) \cap r(g) \not= \emptyset \}} \gamma_\ell
    \sum_{m \in L_{\ell,2}} \lambda_{\ell,m} x^*_{\ell,m} (g) \biggr| \label{eqn:M} \\
    & \quad + \frac{47}{m_{j_0}} \biggl(\sum_{\{\ell \in B :r(x^*_\ell) \cap r(g) \not= \emptyset \}}
    \gamma_\ell^2 \biggr)^{\frac{1}{2}} \nonumber
    + 6 \biggl(\sum_{\{\ell \in E : r(x^*_\ell) \cap r(g) \not= \emptyset \}} \gamma_\ell^2 \biggr)^{\frac{1}{2}}. \nonumber
\end{align}

\n For the second inequality we applied Lemma \ref{lem:geq}.  Now
since $(x^*_{\ell,m})_m$ is $S_{f_{j_0}}$ admissible and
$(x_\ell^*)_\ell$ is $S_{n_i}$ admissible we can use the convolution
property of Schreier families to conclude that
$((x^*_{\ell,m})_{\ell \in S, m \in L_\ell}$ is $S_{f_{j_0}+n_i}$
and hence $S_{2f_{j_0}}$ admissible. Thus for $\xi = 2f_{j_0}
<n_{j_0}$ apply Lemma \ref{lem:geq} to the second term of the sum to
obtain,

\begin{align}
&\biggl| \sum_{\{\ell \in S : r(x^*_\ell) \cap r(g) \not= \emptyset \}} \gamma_\ell \sum_{m \in L_{\ell,2}} \lambda_{\ell,m} x^*_{\ell,m} (g) \biggr| \nonumber\\
    &  \leq  \frac{47}{m_{j_0}} \biggl(\sum_{\{\ell \in S : r(x^*_\ell) \cap r(g) \not= \emptyset \} } \gamma^2_\ell \sum_{\substack{m \in L_{\ell,2} \\
    w(x^*_{\ell,m})>m_{j_0}}} \lambda^2_{\ell,m} \biggr)^\frac{1}{2} +
    6 \biggl(\sum_{\{\ell \in S : r(x^*_\ell) \cap r(g) \not= \emptyset \} } \gamma^2_\ell \sum_{\substack{m \in L_{\ell,2} \nonumber\\
    w(x^*_{\ell,m})=m_{j_0}}} \lambda^2_{\ell,m} \biggr)^\frac{1}{2} \nonumber \\
    & \leq  \frac{47}{m_{j_0}} \biggl(\sum_{\{\ell \in S : r(x^*_\ell) \cap r(g) \not= \emptyset \} } \gamma^2_\ell
    \biggl( \frac{1}{w(x^*_\ell)} \biggr)^2 \biggr)^\frac{1}{2} +  6 \biggl(\sum_{\{\ell \in S : r(x^*_\ell) \cap r(g) \not= \emptyset \} }
    \gamma^2_\ell \biggl( \frac{1}{w(x^*_\ell)} \biggr)^2 \biggr)^\frac{1}{2} \nonumber ~~ \mbox{(by (\ref{eqn:K}))} \\
    & \leq  \frac{47}{m_{j_0}\min_\ell \{w(x^*_\ell )\}} \biggl(\sum_{\{\ell \in S : r(x^*_\ell) \cap r(g) \not= \emptyset \} } \gamma^2_\ell
     \biggr)^\frac{1}{2} +  \frac{6}{\min_\ell \{w(x^*_\ell )\}} \biggl(\sum_{\{\ell \in S : r(x^*_\ell) \cap r(g) \not= \emptyset \} }
    \gamma^2_\ell \biggr)^\frac{1}{2} \label{eqn:N}  \\
    & \leq \frac{53}{m_e} \biggl(\sum_{\{\ell \in S : r(x^*_\ell) \cap r(g) \not= \emptyset \} }
    \gamma^2_\ell \biggr)^\frac{1}{2}. \nonumber
\end{align}

Now separate the first term of (\ref{eqn:M}) into two terms. Then
recall $g=y/\|y\|$, $\|y\|\geq 1/m_{j_0}$, that
$((x^*_{\ell,m})_{\ell \in S, m \in L_\ell}$ is $S_{2f_{j_0}}$
admissible and apply Lemma \ref{lem:el} for ``$\xi$''$= 2f_{j_0}<
n_{j_0}$,

$$\mbox{``} \sum_m \gamma_m x^*_m \mbox{''} = \sum_{\{\ell \in S : r(x^*_\ell) \cap
r(g) \not= \emptyset \}} \gamma_\ell \sum_{m \in L_{\ell,1}}
\lambda_{\ell,m} x^*_{\ell,m}$$

\n and ``$\sum_k \beta_k y_k$''$=y$ to conclude that,

\begin{equation}
\begin{split}
\biggl| \sum_{\{\ell \in S : r(x^*_\ell) \cap r(g) \not= \emptyset \}} & \gamma_\ell \sum_{m \not\in L_{\ell,2}} \lambda_{\ell,m} x^*_{\ell,m} (g) \biggr| \\
    & \leq \biggl| \sum_{\{\ell \in S : r(x^*_\ell) \cap r(g) \not= \emptyset \}} \gamma_\ell \sum_{m \in L_{\ell,1}} \lambda_{\ell,m} x^*_{\ell,m} (g) \biggr| \\
    & \quad \quad \quad + m_{j_0} \biggl| \sum_{\{\ell \in S : r(x^*_\ell) \cap r(g) \not= \emptyset \}}
    \gamma_\ell \sum_{m \in L_{\ell,3}} \lambda_{\ell,m} \gamma_{\ell,m} e^*_{p(\ell,m)}
    \sum_k (n_{j_0})_1^R(p_k) y_{j_k}  \biggr| \\
    & \leq 22 m_{j_0} \biggl( \sum_{\{\ell \in S : r(x^*_\ell) \cap r(g) \not= \emptyset \}} \gamma_\ell^2 \sum_{m \in L_{\ell,1}} \lambda_{\ell,m}^2 \biggr)^\frac{1}{2} \\
    & \quad \quad \quad + m_{j_0} \biggl| \sum_{\{\ell \in S : r(x^*_\ell) \cap r(g) \not= \emptyset \}} \gamma_\ell \sum_{m \in L_{\ell,3}} \lambda_{\ell,m} \gamma_{\ell,m} (n_{j_0})_1^R(p_{k(\ell,m )})
    \biggr|,
    \label{eqn:AB}
\end{split}
\end{equation}

\n where for every $\ell \in S$ and $m \in L_{\ell,3}$, $k(\ell,m)$
is the unique integer $k$ (if any) such that
$e^*_{p(\ell,m)}(y_{j_k}) \not=0$.  If no such $k$ exists then the
corresponding term is absent from the second part of the sum.  Now
(\ref{eqn:AB}) continues as follows by applying (\ref{eqn:K}), the
Cauchy-Schwartz inequality and the facts $|\gamma_{\ell,m} | \leq 1$
and $(k(\ell,m))_{\ell \in S, m \in L_{\ell,3}} \in S_{n_i+f_{j_0}}
\subset S_{2f_{j_0}}$ where $2f_{j_0} < n_{j_0}$:

\begin{align}
    & \leq \frac{22 m_{j_0}}{m^2_{j_0}} \biggl( \sum_{\{\ell \in S : r(x^*_\ell) \cap r(g) \not= \emptyset \}}
    \gamma_\ell^2  \biggr)^\frac{1}{2} \nonumber \\
    & \quad \quad + m_{j_0}  \sum_{\{\ell \in S : r(x^*_\ell) \cap r(g) \not= \emptyset \}} |\gamma_\ell |
    \biggl(\sum_{m \in L_{\ell,3}} \lambda^2_{\ell,m} \biggr)^\frac{1}{2} \biggl( \sum_{m \in L_{\ell,3}}
    ((n_{j_0})_1^R(p_{k(\ell,m )}))^2 \biggr)^\frac{1}{2} \nonumber \\
    & \leq \frac{18}{m_{j_0}} \biggl( \sum_{\{\ell \in S : r(x^*_\ell) \cap r(g)
    \not= \emptyset \}} \gamma_\ell^2  \biggr)^\frac{1}{2}
    \nonumber \\
    & \quad \quad + \frac{m_{j_0}}{\min_\ell \{w(x_\ell^*) \}}  \biggl(\sum_{\{\ell \in S : r(x^*_\ell) \cap r(g) \not= \emptyset \}} \gamma^2_\ell \biggr)^\frac{1}{2}
    \biggl( \sum_{\{\ell \in S : r(x^*_\ell) \cap r(g) \not= \emptyset \}} \sum_{m \in L_{\ell,3}}
    ((n_{j_0})_1^R(p_{k(\ell,m )}))^2 \biggr)^\frac{1}{2}     \label{eqn:P} \\
    & \leq  \frac{18}{m_{j_0}} \biggl( \sum_{\{\ell \in S : r(x^*_\ell) \cap r(g) \not= \emptyset \}} \gamma_\ell^2  \biggr)^\frac{1}{2}
    +  \frac{1}{m_{j_0}\min_\ell \{w(x_\ell^*) \}} \biggl( \sum_{\{\ell \in S : r(x^*_\ell) \cap r(g) \not= \emptyset \}} \gamma_\ell^2
    \biggr)^\frac{1}{2} \nonumber \\
    & \leq \frac{23}{m_{j_0}} \biggl( \sum_{\{\ell \in S : r(x^*_\ell)
    \cap r(g) \not= \emptyset \}} \gamma_\ell^2
    \biggr)^\frac{1}{2}. \nonumber
\nonumber
\end{align}

\n The result follows by combining (\ref{eqn:M}),(\ref{eqn:N}),(\ref{eqn:AB}) and (\ref{eqn:P}).

\end{proof}

\begin{cor}
Let $(y_j)_{j\in I}$, $j_0 \in 2\N$ and a normalized
$(1/m^2_{j_0},n_{j_0})$ squared average $g$ of $(y_j)_{j \in I}$
chosen as in Remark \ref{rem:barb}. Then for any $x^* \in \NC$ such
that $x^*(g)
> 1/2$ we have that $w(x^*)=m_{j_0}$. \label{cor:norming}
\end{cor}

\begin{proof}
Let $x^* \in \NC$ such that $x^*(g)>1/2$. Assume
$w(x^*)\not=m_{j_0}$. Apply Lemma \ref{lem:any} for $i=0$, for
(recall that $n_0=0$) the $S_{n_i}$ admissible family
$(x_\ell^*)_\ell$ being the singleton $\{x^*\}$ and $\gamma_1=1$, to
obtain that

\begin{equation}
 x^*(g) < \frac{123}{m_e} \leq \frac{123}{m_1} < \frac{1}{2},
 \label{star}
\end{equation}

\n (where $m_e=\min(m_{j_0},w(x^*))$), since $m_1>246$.
\end{proof}

\section{$X$ is a Hereditarily Indecomposable Banach Space}

We now show that $X$ is HI.  We proceed by fixing $j\in \N$ and by
defining vectors $(g_i)_{i=1}^p$ and $(z_i)_{i=1}^p$, functionals
$(x_i^*)_{i=1}^p \in \NC$, positive integers $(j_i)_{i=1}^p$ and
$R=(t_i) \in \N$ which will be fixed throughout the section and
shall be referred to in the results of the section. By using
standard arguments we can assume that any two subspaces, in our case
with trivial intersection, are spanned by normalized block bases of
$(e_n)$. Let $(u_n)$ and $(v_n)$ be two such normalized block bases
of $(e_n)$ and fix $j \in \mathbb{N}$. Set $P=(p_n)$ and $Q=(q_n)$
where $p_n=\min \supp u_n$ and $q_n=\min \supp v_n$ for all $n \in
\N$. By passing to subsequences of $(p_n)$ and $(q_n)$ and
relabeling, assume by (\ref{eqn:AA}) that if $R \in [P\cup Q]$ then
for $\xi < n_{2j+1}$, $\sup \{ (\sum_{k\in F} ((n_{2j+1}))_1^R(k))^2
)^\frac{1}{2} : F \in S_{\xi} \} < 1/m^2_{2j+1}$. By Lemma
\ref{lem:sn} let $(y_j)_{j\in 2\N -1}$ (resp. $(y_j)_{j \in 2 \N}$)
be a block sequence of $(u_n)_{n}$ (resp. $(v_n)_{n}$) such that
$y_j$ is a smoothly normalized $(1/m_{2j},f_{2j} +1)$ squared
average of $(u_n)_{n}$ (resp. $(v_n)_{n}$). Apply Lemma
\ref{lem:lwl} to $(y_j)_{j\in 2\N -1}$ and $(y_j)_{j \in 2\N }$ to
obtain $I_1 \in [2\N -1]$ and $I_2 \in [ 2\N ]$ such that $(y_j)_{j
\in I_1}$ and $(y_j)_{j \in I_2}$ satisfy the statement of Lemma
\ref{lem:lwl}.  For $j_1 \in \N$, $2j_1>2j+1$ let $g_1$ be a
normalized $(1/m^2_{2j_1},n_{2j_1})$ squared average of $(y_j)_{j
\in I_1}$. Let $x^*_1 \in \NC$ with $x^*_1(g_1)>1/2$ and $r(x^*_1)
\subset r(g_1)$.  By Corollary \ref{cor:norming} we have that
$w(x^*_1)=m_{2j_1}$.  Let $m_{2j_2}=\sigma(x^*_1)$.  Let $g_2$ be a
normalized $(1/m^2_{2j_2},n_{2j_2})$ squared average of $(y_j)_{j
\in I_2}$ with $g_1<g_2$.  Let $x^*_2 \in \NC$  with
$x^*_2(g_2)>1/2$ and $r(x^*_2) \subset r(g_2)$.  By Corollary
\ref{cor:norming} we have that $w(x^*_2)=m_{2j_2}$.  Let
$m_{2j_3}=\sigma(x^*_1,x^*_2)$. Let $g_3$ be a normalized
$(1/m^2_{2j_3},n_{2j_3})$ squared average of $(y_j)_{j \in I_1}$
with $g_1<g_2<g_3$. Let $x^*_3 \in \NC$  with $x^*_3(g_3)>1/2$ and
$r(x^*_3) \subset r(g_3)$. By Corollary \ref{cor:norming} we have
that $w(x^*_3)=m_{2j_3}$. Continue similarly to obtain
$(g_i)_{i=1}^p$ and $(x^*_i)^p_{i=1} \subset \NC$ such that:
\begin{itemize}
\item[(a)] $g_i$ is a normalized $( 1/
m^2_{2j_i},n_{2j_i} )$ squared average of $(y_j)_{j \in I_1}$ for
$i$ odd (resp. $(y_j)_{j\in I_2}$ for $i$ even).
\item[(b)] $w(x^*_i)=m_{2j_i}, r( x^*_i) \subset r(g_i)$ and
$x^*_i(g_i) > 1/2$.
\item[(c)] $\sigma(x^*_1, \cdots , x^*_{i-1})=w(x^*_i)$ for all $2\leq
i\leq p$.
\item[(d)] $\{g_i :i \leq p\}$ is maximally $S_{n_{2j+1}}$- admissible.
\end{itemize}

\n Let $z_i=g_i/(x^*_{i}(g_i))$.  Let $t_i= \min \supp z_i$ and
$R=(t_i)_i$. The fact that $X$ is HI will follow from the next
proposition.

\begin{prop}
For all $x^* \in \NC$ there exist $J_1 < \cdots < J_s \subset \{1,
\cdots ,p\}$ such that,
\begin{enumerate}
\item $\{z_{min J_m} : m\leq s\}\in S_0 + 3 S_{f_{2j+1}}$ (i.e. it can be written
as a union of four sets: one is a singleton and three belong to $S_{f_{2j+1}}$).
\item There exists $(b_m)_{m =1}^s \subset \mathbb{R}$ such that

$$ \left|x^*\biggl(\sum_{i \in J_m} (-1)^i (n_{2j+1})_1^R(t_i)
z_i \biggr)\right| \leq  (n_{2j+1})_1^R(t_{\min J_m})b_m
~\mbox{and}~ \biggl(\sum_m b^2_m \biggr)^\frac{1}{2} \leq 6.$$
\item
$$ \left|x^*\biggl(\sum_{i \in \{1,\ldots, p\}\setminus \cup_{m=1}^s J_m } (-1)^i (n_{2j+1})_1^R(t_i)
z_i \biggr)\right| \leq \frac{505}{m^2_{2j+1}}.$$
\end{enumerate}
\label{prop:conds}
\end{prop}

Before presenting the proof of this proposition we show that it
implies that $X$ is HI.  First we find a lower estimate for
$\|\sum_{i=1}^p (n_{2j+1})_1^R(t_i)z_i \|$.

\begin{equation*}
\begin{split}
\|\sum_{i} (n_{2j+1})_1^R(t_i)z_i \| & \geq
\frac{1}{m_{2j+1}} \sum_k (n_{2j+1})_1^R(t_k) x^*_k \biggl( \sum_{i} (n_{2j+1})_1^R(t_i)z_i \biggr) \\
    & = \frac{1}{m_{2j+1}} \sum_k ((n_{2j+1})_1^R(t_k))^2 = \frac{1}{m_{2j+1}}.
\end{split}
\end{equation*}

\n Now we find an upper estimate for $\|\sum_{i=1}^p
(-1)^i(n_{2j+1})_1^R(t_i)z_i \|$. Let $x^* \in \NC$ and find $J_1 <
J_2 < \cdots < J_s$ to satisfy Proposition \ref{prop:conds}.

\begin{equation*}
\begin{split}
x^*\biggl( &\sum_{i} (-1)^i (n_{2j+1})_1^R(t_i)z_i\biggr) \\
    & \leq \biggl| x^*\biggl(\sum_{i \in \cup_{m=1}^s J_m} (-1)^i (n_{2j+1})_1^R(t_i)z_i\biggr)\biggr| +
    \biggl| x^*\biggl(\sum_{i \not\in \cup_{m=1}^s J_m} (-1)^i (n_{2j+1})_1^R(t_i)z_i\biggr) \biggr| \\
    & \leq \sum_{m=1}^s \biggl| x^* \biggl(\sum_{i \in J_m} (-1)^i (n_{2j+1})_1^R(t_i)z_i\biggr)\biggr|
    + \frac{505}{m^2_{2j+1}} \\
    & \leq \sum_{m=1}^s (n_{2j+1})_1^R(t_{\min J_m}) b_m + \frac{505}{m^2_{2j+1}}  \\
    & \leq \biggl(\sum_{m=1}^s ((n_{2j+1})_1^R(t_{\min J_m}))^2 \biggr)^\frac{1}{2}
    \biggl(\sum_{m=1}^s b_m^2 \biggr)^\frac{1}{2}
    + \frac{505}{m^2_{2j+1}}  \\
    & \leq \biggl(4 \biggl( \frac{1}{m^2_{2j+1}} \biggr)^2
    \biggr)^\frac{1}{2} 6 + \frac{505}{m^2_{2j+1}} =
    \frac{517}{m^2_{2j+1}}
\end{split}
\end{equation*}

\n where the numbers ``4'' and ``6'' after the last inequality are justified
by parts (1) and (2) of Proposition \ref{prop:conds} respectively. 
Combining the two estimates we have that,

$$\frac{517}{m_{2j+1}}\biggl\|\sum_{i=1}^p (n_{2j+1})_1^R(t_i)z_i \biggr\| \geq \frac{517}{m^2_{2j+1}}
\geq \biggl\|\sum_{i=1}^p (-1)^i (n_{2j+1})_1^R(t_i)z_i \biggr\|, $$

\n for any $j$ and thus $X$ is HI.

The remainder of the paper will be devoted to proving Proposition
\ref{prop:conds}. The following three lemmas will be needed in the
proof.

\begin{lem}
Let $x^*\in \NC$ with $w(x^*) > m_{2j+1}$.  Let $V \subset \{1,
\ldots ,p \}$ such that $w(x^*) \not\in \{m_{2j_i} : i \in V \}$.
Then

$$x^* \biggl(\sum_{i \in V} \alpha_i z_i \biggr) < \frac{496}{m^2_{2j+1}}
\biggl(\sum_{\{i \in V : r(x^*)\cap r(z_i)\not= \emptyset
\}}\alpha^2_i \biggr)^\frac{1}{2}.$$

 \label{lem:bw}
\end{lem}

\begin{proof}
Let $x^*=\frac{1}{w(x^*)} \sum_\ell \gamma_\ell y^*_\ell$ where
$\sum_\ell \gamma_\ell^2 \leq 1$ and $(y_\ell^*)_\ell$ is
appropriately admissible. Define,

$$Q(1)=\{i \in V :\mbox{there is exactly one}~ \ell ~\mbox{such that}~ r(y_\ell^*)\cap r(z_i)\not= \emptyset\},$$
$$Q(2)=\{i \in V :\mbox{there is at least two $\ell $'s such that}~
r(y_\ell^*)\cap r(z_i)\not= \emptyset\}.$$

\n For $i$'s in $Q(1)$,

\begin{equation}
\begin{split}
x^* \biggl(\sum_{ i \in Q(1) } \alpha_i z_i \biggr) & = \frac{1}{w(x^*)} \sum_\ell \gamma_\ell
y^*_\ell \biggl( \sum_{\{ i \in Q(1) : r(y_\ell^*)\cap r(z_i)\not= \emptyset \} } \alpha_i z_i \biggr) \\
    & \leq \frac{1}{w(x^*)} \sum_\ell |\gamma_\ell | 2 \sqrt{3}
    \biggl(\sum_{\{ i \in Q(1) : r(y_\ell^*)\cap r(z_i)\not= \emptyset \}} \alpha_i^2 \biggr)^\frac{1}{2} \\
    & \leq \frac{2 \sqrt{3}}{w(x^*)}
    \biggl(\sum_{ \{ i : r(x^*)\cap r(z_i)\not= \emptyset\} } \alpha^2_i \biggr)^\frac{1}{2}
    \leq \frac{4}{m_{2j+1}^2} \biggl( \sum_{\{ i : r(x^*)\cap r(z_i)\not= \emptyset \} } \alpha^2_i
    \biggr)^\frac{1}{2}, \label{eqn:Q}
\end{split}
\end{equation}

\n where for the first inequality we used Proposition \ref{prop:ul2}
and the fact that $x^*_i(g_i)>1/2$, for the second inequality we
applied the Cauchy-Schwartz inequality and for the last inequality
we used the fact that if $2j+1<i$ then $m^2_{2j+1} \leq m_i$. For the
$Q(2)$ case, if $w(x^*) \not= m_{2j_i}$ for all $i \leq p$ notice
that

\begin{equation}
\begin{split}
x^*  & \biggl(\sum_{\{i \in Q(2):r(x^*)\cap r(z_i)\not= \emptyset\}}
\alpha_i z_i \biggr) \\
    & = \sum_{\{i \in Q(2):r(x^*)\cap r(z_i)\not= \emptyset\}} \alpha_i
    \biggl(\biggl( \sum_{\{\ell : r(y^*_\ell)\cap r(z_i)\not= \emptyset \}}
    \gamma_\ell^2 \biggr)^\frac{1}{2}+ \eta \biggr) \biggl|
    \frac{1}{w(x^*)} \sum_{\{\ell :r(y^*_\ell)\cap r(z_i)\not= \emptyset \}}
    \beta_{i,\ell} y^*_\ell (z_i) \biggr| \label{eqn:R}
\end{split}
\end{equation}

\n where $\beta_{i,\ell} = \gamma_\ell / ((\sum_{\{m : r(y_m^*)\cap
r(z_i)\not= \emptyset \}} \gamma_m^2 )^\frac{1}{2}+\eta_i)$ where
$\eta_i$ is arbitrarily small and such that

$$\biggl(\sum_{\{m : r(y_m^*) \cap r(z_i) \not= \emptyset \}} \gamma_m^2
\biggr)^\frac{1}{2}+\eta_i \in \Q.$$

Now apply Lemma \ref{lem:any} for ``$g$''$=g_i$, ``$j_0$''$=2j_i$,
``$i$''$=0$, ``$\gamma_1$''$=1$ and

\n ``$x_1^*$''$=\frac{1}{w(x^*)}\sum_{\{\ell : r(y_\ell^*)\cap
r(z_i)\not= \emptyset \}} \beta_{i,\ell} y^*_\ell \in \NC$ to
continue (\ref{eqn:R}) as follows:

\begin{equation}
\begin{split}
& \leq \sum_{\{i \in Q(2):r(y_\ell^*)\cap r(z_i)\not= \emptyset\}} |
\alpha_i | \biggl(\biggl(\sum_{\{ \ell :r(x^*)\cap r(z_i)\not=
\emptyset\}}
    \gamma_\ell^2 \biggr)^\frac{1}{2} + \eta_i \biggr) \frac{2(123)}{w(x^*)} \\
    & \leq \frac{2(246)}{m^2_{2j+1}} \biggl(\sum_{\{i \in V :r(x^*)\cap r(z_i)\not= \emptyset\}} \alpha_i^2\biggr)^\frac{1}{2},
\label{eqn:R2}
\end{split}
\end{equation}

\n where the constants $\eta_i$ were forgotten in the last
inequality since they were arbitrarily small and the constant
``123'' that appears in the statement of Lemma \ref{lem:any} is
multiplied by 2 since $z_i=g_i/x^*_i(g_i)$, $x^*_i(g_i)>1/2$ and by
another factor 2 since for each $\ell$ there are at most two values
of $i \in Q(2)$ such that $r(y_\ell^*)\cap r(z_i)\not= \emptyset$.
The result follows by combining (\ref{eqn:Q}), (\ref{eqn:R}) and (\ref{eqn:R2}).
\end{proof}

\begin{lem}
Let $x^* \in \NC$ with $w(x^*)=m_{2j+1}$. Thus
$x^*=\frac{1}{m_{2j+1}} \sum_\ell \gamma_\ell y^*_{\ell}$ where
$\sum_\ell \gamma_\ell^2 \leq 1$, $y^*_\ell \in \NC$, and
$(y_\ell^*)_\ell$ has $S_{n_{2j+1}}$ dependent extension. Let $V
\subset \{1 \leq i \leq p: m_{2j_i} \not= w(y^*_\ell)  ~\mbox{for
all}~ \ell \}$ and $(\alpha_i)_i \in c_{00}$. Then

$$x^* \biggl(\sum_{i\in V} \alpha_i z_i \biggr) < \frac{2}{m^2_{2j+1}}
\biggl(\sum_{\{i \in V : r(x^*) \cap r(z_i) \not= \emptyset \}}
\alpha_i^2 \biggr)^\frac{1}{2}.$$ \label{lem:sw}
\end{lem}

\begin{proof}
For $i \in I$ define $Q(1)$ and $Q(2)$ as in Lemma \ref{lem:bw}. For
$i \in Q(1)$, use Lemma \ref{lem:bw} for ``$x^*$''$=y^*_\ell$ for
each $\ell$ to obtain,

\begin{align}
x^*\biggl(\sum_{i\in Q(1)} \alpha_i z_i \biggr) & = \frac{1}{m_{2j+1}} \sum_\ell \gamma_\ell
y^*_\ell \biggl(\sum_{i\in Q(1)} \alpha_i z_i \biggr) \nonumber  \\
    & \leq \frac{1}{m^2_{2j+1}} \sum_\ell |\gamma_\ell |
    \biggl(\sum_{\{i \in V: r(y_\ell^*) \cap r(z_i) \not= \emptyset\}}\alpha^2_i \biggr)^\frac{1}{2} \frac{496}{m^2_{2j+1}}
    \leq \frac{1}{m^2_{2j+1}}
    \biggl(\sum_{\{i \in V: r(x^*) \cap r(z_i) \not= \emptyset\}} \alpha_i^2
    \biggr)^\frac{1}{2}, \label{eqn:ZZ}
\end{align}

\n where for the last inequality we used the Cauchy Schwartz
inequality and the fact that $m_{2j+1} \geq m_3 \geq m_1^4 \geq
496$. For $i \in Q(2)$,

\begin{equation}
\begin{split}
x^* \biggl(\sum_{i\in Q(2)} \alpha_i z_i \biggr) & =
\frac{1}{m_{2j+1}} \sum_\ell \gamma_\ell y^*_\ell \biggl(\sum_{i\in
Q(2)} \alpha_i z_i \biggr) \leq \frac{1}{m_{2j+1}} \sum_{i \in
Q(2)}| \alpha_i |
    \biggl| \sum_{\ell} \gamma_\ell y^*_\ell (z_i) \biggr| \\
    & \leq \frac{1}{m_{2j+1}} \sum_{i\in Q(2)} | \alpha_i |
    \biggl(\sum_{\{\ell : r(y_\ell^*) \cap r(z_i) \not= \emptyset\}} \gamma^2_i \biggr)^\frac{1}{2}
    \frac{2(123)}{\min_{\ell}\{w(y^*_\ell), m_{2j_i})\}} \\
    & \leq \frac{4(123)}{m^3_{2j+1}} \biggl(\sum_{\{i :r(x^*) \cap r(z_i) \not= \emptyset\}} \alpha_i^2
    \biggr)^\frac{1}{2} \leq \frac{1}{m^2_{2j+1}} \biggl(\sum_{\{i :r(x^*) \cap r(z_i) \not= \emptyset
    \}}
    \alpha_i^2 \biggr)^\frac{1}{2}. \label{eqn:ZZZ}
\end{split}
\end{equation}

\n Recall that $(y_\ell^*)_\ell$ is $S_{n_{2j+1}}$ admissible and
$2j+1 <2j_i$ for all $i$.  Thus the second inequality follows from
applying Lemma \ref{lem:any} for ``$x^*_\ell$''$=y^*_\ell$,
``$i$''$=2j+1$, ``$j_0$''$=2j_i$, ``$g$''$=g_i$ and observing that
$z_i = g_i / x^*_i(g_i)$, $x^*(g_i) >1/2$ and $\{w(y^*_\ell):\ell\}
\cap \{m_{2j_i} : i \in V\} = \emptyset$.  The third inequality
follows from applying the Cauchy-Schwartz inequality and observing
that for every $\ell$ there are at most two values of $i \in Q(2)$
such that $r(y_\ell^*) \cap r(z_i) \not= \emptyset $ and
$\min_{i,\ell} \{w(y^*_{\ell}),m_{2j_i}\} \geq m^2_{2j+1}$.  The
fourth inequality follows since $m_{2j+1} \geq m_3 \geq m_1^4 \geq
492$.  The result follows by combining (\ref{eqn:ZZ}) and
(\ref{eqn:ZZZ}).

\end{proof}

\begin{lem}
For $x^* \in \NC$, $w(x^*)=m_{2j+1}$ there exist $J_1 < J_2 < J_3$
subsets of $\{1, \ldots , p \}$ (some of which are possibly empty)
such that:
\begin{enumerate}
\item For $m\in \{1,2,3\}$,

$$\biggl|x^* \biggl( \sum_{i \in J_m} (-1)^i (n_{2j+1})_1^R (t_i)
z_i \biggr) \biggr| \leq 2 (n_{2j+1})_1^R(t_{min J_m}).$$

\item $x^* = \frac{1}{m_{2j+1}} \sum_{\ell=1}^n \gamma_\ell y_\ell^*$ with
$\sum_\ell \gamma_\ell^2 \leq 1$, $y_\ell^* \in \NC$,
$(y_\ell^*)_\ell$ has a $S_{n_{2j+1}}$-dependent extension and
$\{w(y_\ell^*) : 1 \leq \ell \leq n\} \cap \{ m_{2j_i} : i \in
\{1,\ldots,p\} \setminus \cup_{m=1}^3 J_m \} = \emptyset $.

\item Also,
$$\biggl|x^* \biggl( \sum_{i \in
\{1,\ldots,p\} \setminus \cup_{m=1}^3 J_m} (-1)^i (n_{2j+1})_1^R
(t_i) z_i \biggr) \biggr| \leq \frac{2}{m^2_{2j+1}}.$$
\end{enumerate}
Moreover, for any interval $Q \subset \{1, \ldots, p \}$ there exist
$J_1 < J_2 < J_3$ subsets of $Q$ (some of which are possibly empty)
such that conditions (1), (2) and (3) are satisfied with the
exception that in conditions (2) and (3) the set $\{1, \ldots ,p
\}\setminus \cup_{m=1}^3 J_m $ is replaced by $Q\setminus
\cup_{m=1}^3 J_m $.
 \label{lem:new}
\end{lem}

\begin{proof}
Suppose $x^* \in \NC$ with $w(x^*)=m_{2j+1}$.  Then $x^*=1/m_{2j+1}
\sum_{\ell=1}^n \gamma_\ell y^*_\ell$ with $\sum_\ell \gamma_\ell^2
\leq 1$,  $y_\ell^* \in \NC$ and $(y^*_\ell)_{\ell=1}^n$ has a
$S_{n_{2j+1}}$-dependent extension. Thus there exists $d,L \in \N$
and an $S_{n_{2j+1}}$ admissible family $(\tilde y_\ell^*)_{\ell
=1}^{d+n-1}$ such that $\sigma( \tilde y_1^*, \ldots , \tilde
y^*_\ell) = w(\tilde y^*_{\ell+1}) $ for $1 \leq \ell < d+n-1$ and
$\tilde y^*_{\ell} |_{[L,\infty)} = y^*_{\ell-(d-1)}$ for $\ell = d,
\ldots , d+n-1$.

Recall the definition of $(x^*_k)$ from the beginning of this
section.  By injectivity of $\sigma$, the set $\{k \in \{1, \ldots
,p \} :w(x^*_k) \in \{w(\tilde y_\ell^*) : d \leq i \leq d+n-1\}\}$
is an interval of integers (possibly empty). Let $k_0$ be the
largest integer $k$ such that $w(x_k^*) \in \{w(\tilde y_i^*):d \leq
i \leq d+n-1 \}$ and $k_0=0$ if no such $k$ exists.

If $k_0=0$ (i.e. $w(x^*_k)\not=w(y_\ell^*)$ for all $1 \leq \ell\leq
n,1 \leq k\leq p$) then let $J_1=J_2=J_3= \emptyset$ and conditions
(1) and (2) are trivial.  To verify condition (3) apply Lemma
\ref{lem:sw} for ``$V$''$=\{k: 1 \leq k \leq p\}$ and
``$\alpha_i$''$=(-1)^i(n_{2j+1})_1^R(t_i)$.

If $k_0=1$ (i.e. $w(x^*_1) = w(y_{i_0}^*)$ for some $i_0 \in \{1,
\ldots, n\}$ and $w(x^*_k) \not\in \{w(y_\ell^*):1 \leq \ell \leq
n\}$ for $1 < k \leq p$) then set $J_1=\{1\}$, $J_2=J_3=\emptyset$
and since $\|z_1\|=\|g_1/x^*_1(g_1)\| \leq 2$ it is easy to check
that conditions (1) and (2) hold.  To check condition (3) apply
Lemma \ref{lem:sw} for ``$V$''$= \{2,3,\ldots,p \}$ and
``$\alpha_i$''$=(-1)^i(n_{2j+1})_1^R(t_i)$.

If $k_0 >1$ and $w(x^*_{k_0}) = w(\tilde y_1^*)$ then by the
injectivity of $\sigma$, $w(x^*_k) \not\in \{w(y_\ell^*): 1\leq \ell
\leq n \}$ for $k \in \{1, \ldots, p\}\setminus \{k_0\}$. Thus set
$J_1=\{k_0\}$, $J_2=J_3= \emptyset$ and easily verify conditions (1)
and (2) as above. To check condition (3) apply Lemma \ref{lem:sw}
for ``$V$''$=\{1, \ldots,p\}\setminus \{k_0\}$.

Finally, if $k_0 >1$ and $w(x^*_{k_0})=w(\tilde y_{\ell_0}^*)$ for
some $\ell_0 >1$ by the injectivity of $\sigma$ it must be the case
that $\ell_0=k_0$, $w(\tilde y_i^*) =w(x^*_i)$ for all $i \leq k_0$
and $\tilde y_i^* = x^*_i$ for all $i < k_0$.  Let $J_1=\{d\}$ if $d
\leq k_0$ (else $J_1=J_2=J_3=\emptyset$), $J_2=(d,k_0)\cap \N$ if $d
< k_0$ (else $J_2=J_3=\emptyset$), $J_3=\{k_0\}$ if $d<k_0$. By the
choice of $J_1, J_2, J_3$ we have that $w(x_k^*) \not\in
\{w(y_\ell^*): 1 \leq \ell \leq n \}$ for $k \in \{1, \ldots ,p\}
\setminus \cup_{m=1}^3 J_m$ and if $J_2 \not=\emptyset$ then
$x^*_i=\tilde y^*_i = y^*_{i-d+1}$ for $i \in J_2$.  Apply Lemma
\ref{lem:sw} for ``$V$''$=\{1, \ldots ,p \}\setminus \cup_{m=1}^3
J_m$ and ``$\alpha_i$''$=(-1)^i(n_{2j+1})_1^R(t_i)$ to satisfy
conditions (2) and (3).

To verify condition (1) for $J_1$ and $J_3$ (if they are non-empty),
since they are singletons, simply observe that $|x^*(z_i)|\leq
\|z_i\| = \|g_i/x^*_i(g_i)\| \leq 2$. To verify conditions (1) for
$J_2$ (if it is non-empty),

\begin{equation}
\begin{split}
\biggl| x^* \biggl( \sum_{i \in J_2} (-1)^i(n_{2j+1})_1^R (t_i) z_i
\biggr) \biggr| & = \frac{1}{m_{2j+1}} \biggl| \sum_{i \in J_2}
\gamma_{i-(d-1)} y^*_{i-(d-1)} \biggl( \sum_{i \in J_2} (-1)^i
(n_{2j+1})_1^R (t_i) z_i \biggr)  \biggr| \\
    & =  \frac{1}{m_{2j+1}} \biggl| \sum_{i \in J_2} \gamma_{i-(d-1)}
    x_i^*
    \biggl( \sum_{i \in J_2} (-1)^i (n_{2j+1})_1^R (t_i) z_i \biggr)  \biggr| \\
    & = \frac{1}{m_{2j+1}} \biggl| \sum_{i \in J_2} \gamma_{i-(d-1)} (-1)^i (n_{2j+1})_1^R (t_i) \biggr| \\
    & \leq \frac{\gamma_{\min J_2-(d-1)}}{m_{2j+1}} (n_{2j+1})_1^R(t_{\min J_2}) \leq (n_{2j+1})_1^R(t_{\min J_2}),
    \label{eqn:Z1}
\end{split}
\end{equation}

\n where the first two equalities follow from the fact that $x_i^*=
y^*_{i-(d-1)}$ for $i \in J_2$; the third equality follows from the
fact that $x^*_i(z_i)= x^*_i(g_i/x^*_i(g_i))=1$; the first
inequality follows from the fact that $((n_{2j+1})_1^R(t_i))_i$ and
$(\gamma_i)_i$ are both non-increasing sequences of non-negative
numbers and $J_2$ is an interval.

\n The proof of the moreover part is identical to the above with the
only exception that the sets $J_1 , J_2 , J_3$ chosen above are
replaced by $Q \cap J_1, Q \cap J_2, Q \cap J_3$.  Notice it was
important in the proof of (\ref{eqn:Z1}) that $J_2$ was an interval.
This remains true if $J_2$ is replaced by $J_2 \cap Q$ since $Q$ was
assumed to be an interval.

\end{proof}

\proof[Proof of Proposition \ref{prop:conds}] Suppose $w(x^*) >
m_{2j+1}$, apply Lemma \ref{lem:bw} for

``$\alpha_i$''$=(-1)^i (n_{2j+1})_1^R(t_i)$.  The conclusion 
of Propostion \ref{prop:conds} is satisfied with $s=1$ and $J_1=\{q\}$ if $w(x^*)=m_{2j_q}$ or $J_1 =
\emptyset$ if $w(x^*) \not\in \{m_{2j_i} : 1 \leq i \leq p \}$.

If $w(x^*)=m_{2j+1}$ the proposition follows directly from Lemma \ref{lem:new}
with $b_1=b_2=b_3=2$.

Assume $w(x^*)<m_{2j+1}$.  Write $x^*=\sum_{\ell \in L} \gamma_\ell
y^*_\ell$ where $(y_\ell^*)_\ell$, $(\lambda_\ell)_\ell$ and
$L=L_1\cup L_2 \cup L_3$ satisfy the conclusions of Lemma
\ref{lem:dl} for ``$j$''$=2j+1$.  Let $L_2=L_2^\prime \cup L_2^{\prime\prime}$ where
$L_2^{\prime}=\{\ell \in L_2 : w(y^*_\ell)=m_{2j+1} \}$ and
$L_2^{\prime\prime}=\{\ell \in L_2 : w(y^*_\ell)>m_{2j+1} \}$.
Define,

$$Q(1) = \{1 \leq i \leq p: \mbox{there is exactly one }~ \ell \in L ~\mbox{with}~~
r(y^*_\ell)\cap r(z_i)\not= \emptyset \},$$

$$Q(2)=\{1 \leq i \leq p: \mbox{there are at least two $\ell $'s in $L$ with}~~
r(y^*_\ell)\cap r(z_i)\not= \emptyset\}.$$

\n For $i$'s in $Q(2)$,

\begin{align}
\biggl|\biggl(\sum_{\ell \in L} \lambda_\ell y^*_\ell \biggr)
\biggl(\sum_{i\in Q(2)} & (-1)^i (n_{2j+1})_1^R (t_i) z_i\biggr)
\biggr|  \leq
\sum_{i\in Q(2)} (n_{2j+1})_1^R (t_i) \biggl| \sum_{\ell \in L} \lambda_\ell y^*_\ell(z_i) \biggr| \nonumber \\
    & \leq 2 \sum_{i\in Q(2)} (n_{2j+1})_1^R (t_i) \biggl| \sum_{\ell \in L} \lambda_\ell
    y^*_\ell(g_i) \biggr| \label{eqn:q2}\\
    & \leq 2(7)  \sum_{i\in Q(2)}  (n_{2j+1})_1^R (t_i) \biggl(\sum_{\{
    \ell \in L :r(y^*_\ell) \cap r(z_i) \not= \emptyset  \}}
    \lambda^2_\ell \biggr)^{\frac{1}{2}}  \nonumber \\
    & \leq 14 \biggl( \sum_{i \in Q(2)} ((n_{2j+1})_1^R(t_i))^2
    \biggr)^\frac{1}{2} \biggl( \sum_{i \in Q(2)} \sum_{\{ \ell \in
    L : r(y_\ell^*) \cap r(z_i) \not= \emptyset \}} \lambda_\ell^2
    \biggr)^\frac{1}{2} \nonumber \\
    & \leq 14 \frac{2}{m_{2j+1}^2} \frac{2}{w(x^*)} \leq
    \frac{1}{m_{2j+1}^2} \nonumber,
\label{eqn:q2}
\end{align}

\n where for the second inequality we used the definition of $z_i$.
For the third inequality we applied Lemma  \ref{lem:any} for
``$j_0$''$=2j_i$ and ``$g$''$=g_i$ for each $i$ using the fact that
$(y_\ell^*)_\ell$ is $S_{f_{2j+1}}$ admissible, $f_{2j+1} <
n_{2j+1}$ and $2j+1 < 2j_i$ and noticing that the right hand side of
(\ref{eqn:star}) is at most equal to $7(\sum_{\{\ell:
r(x^*_{\ell})\cap r(g) \not= \emptyset \}} \gamma_\ell^2
)^\frac{1}{2}$ since $103/m_e \leq 1$. For the fourth inequality we
used the Cauchy-Schwartz inequality. For the fifth inequality we
used that $(t_i)_{i \in Q(2)}$ is $2S_{f_{2j+1}}$ is admissible,
condition (2) of Lemma \ref{lem:dl}, and the fact that for every
$\ell \in L$ there are at most two values of $i \in Q(2)$ such that
$r(y^*_\ell) \cap r(z_i) \not= \emptyset$. Now for $\ell$'s in $L_1$
and $i$'s in $Q(1)$,

\begin{equation}
\begin{split}
\biggl| \biggl(\sum_{\ell \in L_1} \lambda_\ell y^*_\ell \biggr)&
\biggl( \sum_{\substack{i\in Q(1) \\ r(y^*_\ell) \cap r(z_i) \not=
\emptyset }} (-1)^i (n_{2j+1})_1^R (t_i) z_i  \biggr) \biggr| \\
& \leq 2 \biggl| \biggl(\sum_{\ell \in L_1} \lambda_\ell y^*_\ell
\biggr) \biggl(\sum_{\substack{i\in Q(1)\\ r(y^*_\ell) \cap
r(z_i)\not=
\emptyset } } (-1)^i (n_{2j+1})_1^R (t_i) g_i  \biggr) \biggr| \\
    & \leq 2 \sqrt{3} \sum_{\ell \in L_1}
    |\lambda_\ell| \biggl(\sum_{\substack{ i \in Q(1)\\ r(y^*_\ell) \cap r(z_i)\not= \emptyset  }}
    ((n_{2j+1})_1^R (t_i))^2 \biggr)^\frac{1}{2}
    \leq \frac{4}{m^2_{2j+1}},
\label{eqn:l1}
\end{split}
\end{equation}

\n where Proposition \ref{prop:ul2} was used in the second
inequality and for the last inequality we used the Cauchy-Schwartz
inequality and condition (2) of Lemma \ref{lem:dl}.  For $\ell \in
L_3$ there is at most one value of $i$ (call it $i(\ell)$) such that
$r(y^*_\ell) \cap r(z_i)\not= \emptyset $. Thus

\begin{equation}
\begin{split}
\biggl| \biggl( \sum_{\ell \in L_3} \lambda_\ell y^*_\ell \biggr)
\biggl(\sum_{\{i\in Q(1): r(y^*_\ell) \cap r(z_i) \not= \emptyset
\}} (-1)^i (n_{2j+1})_1^R (t_i) z_i \biggr) \biggr| &  \leq
\sum_{\ell \in L_3} | \lambda_\ell |
    (n_{2j+1})_1^R(t_{i(\ell)}) 2 | \\
    & \leq \frac{2}{m^2_{2j+1}},
\label{eqn:l3}
\end{split}
\end{equation}

\n where for the first inequality we used that $|y^*_\ell(z_i)| \leq
\|z_i\| \leq 2$ and for the second inequality we applied the
Cauchy-Schwartz inequality and used that $(t_{i(\ell)})_{\ell \in
L}$ is $S_{f_{2j+1}}$ admissible and $f_{2j+1} <n_{2j+1}$. For
$\ell'' \in L_2''$, set $J_{\ell''}=\emptyset$ if $w(y^*_{\ell''})
\not\in \{m_{2j_i} : 1 \leq i \leq p\}$ and $J_{\ell''}= \{q\}$ if
$w(y_{\ell''}^*)=m_{2j_q}$ for some $q \in \{1, \ldots , p\}$ and
$r(y_{\ell''}^*) \cap r(z_q) \not= \emptyset$. Notice that for
$\ell'' \in L''_2$ if $J_{\ell''}=\{q\} \not= \emptyset$ then we
have,

\begin{equation}
\begin{split}
\biggl|x^*\biggl( \sum_{i \in J_{\ell''}} (-1)^i (n_{2j+1})_1^R(t_i) z_i \biggr)\biggr|
& \leq (n_{2j+1})_1^R(t_q) |x^*(z_q)| \\
    & = (n_{2j+1})_1^R(t_q) |\lambda_{\ell''} y^*_{\ell''}(z_q)| \\
    & \leq (n_{2j+1})_1^R(t_q) |\lambda_{\ell''}| 2 \quad (\text{since}~\|z_q\|\leq 2). \label{eqn:bs1}
\end{split}
\end{equation}

\n For $\ell' \in L_2'$ apply the moreover part of Lemma
\ref{lem:new} for ``$x^*$''$=y^*_{\ell'}$ and ``$Q$''$=Q_{\ell'}=\{i
\in Q(1) : r(z_i) \cap r(y^*_{\ell'}) \not= \emptyset \}$ to obtain
intervals $J_{\ell',1} < J_{\ell',2} < J_{\ell',3} \subset
Q_{\ell'}$ (possibly empty) such that

\begin{equation} \label{eqn:S}
 \biggl| y_{\ell'}^*\biggl( \sum_{i \in J_{\ell',m}} (-1)^i (n_{2j+1})_1^R (t_i) z_i \biggr) \biggr|
\leq 2 (n_{2j+1})_1^R (t_{\min J_{\ell',m}}), \quad m \in \{1,2,3\}
\end{equation}

\n and

\begin{equation}
\begin{split}
& y^*_{\ell'}=\frac{1}{m_{2j+1}} \sum_k \gamma_{\ell',k}
y^*_{\ell',k}~ \mbox{with}~ \sum_k \gamma_{\ell',k}^2 \leq 1,
y^*_{\ell',k} \in \NC, (y^*_{\ell',k})_k ~\mbox{is}~ S_{n_{2j+1}}~~
\mbox{admissible and} \\
& \{w(y^*_{\ell',k}): k \} \cap \{ m_{2j_i} :i \in Q_{\ell'}
\setminus \cup_{m=1}^3 J_{\ell',m} \} = \emptyset. \label{eqn:T}
\end{split}
\end{equation}

\n Thus for $\ell' \in L_2'$ and $m \in \{1,2,3\}$ apply the fact
$J_{\ell',m} \subset Q_{\ell'} := \{i \in Q(1) :r(z_i) \cap
r(y^*_{\ell'}) \not= \emptyset \}$ and (\ref{eqn:S}) to obtain

\begin{equation}
\begin{split}
\biggl|x^*\biggl( \sum_{i \in J_{\ell',m}} (-1)^i (n_{2j+1})_1^R(t_i) z_i \biggr)\biggr|
& \leq |\lambda_{\ell'}| \biggl| y^*_{\ell'}
\biggl( \sum_{i \in J_{\ell,m}} (-1)^i (n_{2j+1})_1^R(t_i) z_i \biggr)\biggr|  \\
    & \leq (n_{2j+1})_1^R(t_{min J_{\ell',m}}) |\lambda_{\ell'}| 2. \label{eqn:bs2}
\end{split}
\end{equation}

\n Thus for $\ell'' \in L_2''$ if $J_{\ell''}\not=\emptyset$, set
$b_{\ell''} =2 \lambda_{\ell''}$ and for $\ell' \in L_2'$ and $m \in\{1,2,3\}$ if
$J_{\ell',m} \not= \emptyset$ set $b_{\ell',m}=2 \lambda_{\ell'}$.  Hence,

\begin{align}
\sum_{\{\ell'' \in L_2'' :J_{\ell''} \not= \emptyset\}} b^2_{\ell''}
+ & \sum_{\{\ell' \in L_2' , m \in \{1,2,3\} : J_{\ell',m}
\not=\emptyset\}} b^2_{\ell',m} \leq \sum_{\{\ell'' \in L_2''
:J_{\ell''} \not=\emptyset\}} 4 \lambda_{\ell''}^2 +
\sum_{\{\ell' \in L_2' , m \in \{1,2,3\} : J_{\ell',m} \not=\emptyset\}} 4 \lambda_{\ell'}^2 \nonumber \\
    &  \leq 16 \sum_{\ell \in L} \gamma^2_\ell \leq \frac{16}{w(x^*)^2} \leq 1 \label{eqn:bs},
\end{align}

\n where the third inequality was obtained from condition (2) of
Lemma \ref{lem:dl}. Index the $b_{\ell''}$ (for $\ell'' \in L_2''$ with
$J_{\ell''} \not= \emptyset$) and $b_{\ell',m}$ (for $\ell' \in L_2'$ and 
$m \in \{1,2,3 \}$ with $J_{\ell',m} \not= \emptyset$) as $(b_m)_{m=1}^s$
the non-empty $J_{\ell''}$ ($\ell''\in
L_2''$) and $J_{\ell',m}$ ($\ell'\in L_2', m \in \{1,2,3\}$) as
$(J_m)_{m=1}^s$ and we can see that (\ref{eqn:bs1}), (\ref{eqn:bs2}) and (\ref{eqn:bs})
imply that condition (2) of Proposition \ref{prop:conds} is satisfied.
Condition (1) of Proposition \ref{prop:conds} follows from the facts that $(y^*_\ell)_{\ell \in L}$
is $S_{f_{2j+1}}$ admissible, for every $\ell'' \in L_2''$ for which
$J_{\ell''}=\{q\} \not= \emptyset$ we have $r(z_q) \cap
r(y_{\ell'}^*) \not= \emptyset$ and for every $\ell' \in L_2'$ and
$m \in \{1,2,3\}$ such that $J_{\ell',m} \not= \emptyset$ we have
$J_{\ell',m} \subset \{i \in Q(1): r(z_i) \cap r(y^*_{\ell'}) \not=
\emptyset\}$ (for $m=1,2,3$). The next Lemma \ref{lem:final} implies
that

\begin{equation}
\begin{split}
& \biggl| \biggr( \sum_{\ell \in L_2} \lambda_\ell y^*_\ell \biggr)
\biggl( \sum_{i \not\in \cup_{m=1}^s J_m} (-1)^i (n_{2j+1})_1^R
(t_i) (z_i) \biggr)\biggr| \leq \frac{498}{m^2_{2j+1}}.
\label{eqn:last}
\end{split}
\end{equation}

\n Indeed apply Lemma \ref{lem:final} for ``$E$''$=L'_2$, $B=L_2''$,
``$V$''$=Q(1) \setminus \cup_{m=1}^s J_m$, ``$\alpha_i$''$=
(-1)^i(n_{2j+1})_1^R(t_i)$. Note that condition (1) of Lemma
\ref{lem:final} is satisfied by (\ref{eqn:T}).  Condition (2) of
Lemma \ref{lem:final} is satisfied by the definitions of $L_2''$ and
$J_{\ell''}$ for $\ell'' \in L_2''$.  Condition (3) of Lemma
\ref{lem:final} is satisfied since $V \subset Q(1)$. Thus equations
(\ref{eqn:q2}), (\ref{eqn:l1}),(\ref{eqn:l3}) and (\ref{eqn:last})
imply condition (3) and finish the proof of Proposition
\ref{prop:conds}.

\begin{flushright}
$Q.E.D.$
\end{flushright}

\begin{lem}
Let $E$ and $B$ be two finite index sets, $V \subset \{1, \ldots ,
p\}$, $(y_\ell^*)_{\ell \in E\cup B} \subset \NC$, $(\lambda_\ell)_{\ell \in
E \cup B}$ be scalars with $\sum_{\ell \in E \cup B} \lambda_\ell^2
\leq 1$ and $(\alpha_i)_{i\in V}$ be scalars with $\sum_{i \in V}
\alpha_i^2 \leq 1$ such that,
\begin{enumerate}
\item For every $\ell \in E$, $y_\ell^*$ can be written as

$$y^*_\ell = \frac{1}{m_{2j+1}} \sum_{\ell =1}^{r_\ell} \gamma_{\ell,s} y^*_{\ell,s}$$

\n with $y^*_{\ell,s} \in \NC$, $(y^*_{\ell,s})_{\ell =1}^{r_\ell}$
is $S_{n_{2j+1}}$-admissible and $\{w(y^*_{\ell,s}): 1 \leq s \leq
r_\ell \} \cap \{ m_{2j_i} :i \in V\} = \emptyset $.
\item For every $\ell \in B$, $w(y_\ell^*) > m_{2j+1}$ and $w(y_\ell^*) \not\in \{m_{2j_i} : i \in
V\}$.
\item For every $k \in V$ there is a unique $\ell \in E \cup B$ such that $r(y_\ell^*) \cap r(z_i) \not= \emptyset $.
\end{enumerate}
Then,

$$\biggl| \biggl( \sum_{\ell \in E\cup B} \lambda_\ell y^*_\ell \biggr) \biggl( \sum_{i \in V}
\alpha_i z_i \biggr) \biggr| \leq \frac{498}{m^2_{2j+1}}.$$
\label{lem:final}
\end{lem}

\begin{proof}
Let $E$, $B$, $(\lambda_\ell)_{\ell \in E \cup B}$,
$(y^*_\ell)_{\ell \in E\cup B}$ and $V$ be given as in the
hypothesis.  Using the triangle inequality we separate into the
cases $\ell \in B$ and $\ell \in E$.  For $\ell$'s in $B$,

\begin{equation*}
\begin{split}
\biggl| \biggl( \sum_{\ell \in B} \lambda_\ell y^*_\ell \biggr)
\biggl( \sum_{\{i \in V : r(y_\ell^*) \cap r(z_i) \not= \emptyset
\}} \alpha_i z_i \biggr) \biggr| & \leq \sum_{\ell\in B}
|\lambda_\ell | \frac{496}{m^2_{2j+1}}
\biggl(\sum_{\{i \in V : r(y_\ell^*) \cap r(z_i) \not= \emptyset\}} \alpha_i^2 \biggr)^\frac{1}{2} \\
    & \leq \frac{496}{m^2_{2j+1}},
\end{split}
\end{equation*}

\n where for the first inequality we used Lemma \ref{lem:bw} for
``$x^*$''$=y^*_\ell$, noting that its hypothesis is satisfied by our
assumption (2) and for the second inequality we used the
Cauchy-Schwartz inequality.  For $ \ell$'s in $E$,

\begin{equation*}
\begin{split}
\biggl| \biggl( \sum_{\ell \in E} \lambda_\ell y^*_\ell \biggr)
\biggl( \sum_{\{i \in V : r(y_\ell^*) \cap r(z_i) \not= \emptyset
\}} \alpha_i z_i \biggr) \biggr| & \leq \sum_{\ell\in E} |
\lambda_\ell | \frac{2}{m^2_{2j+1}}
\biggl(\sum_{\{i \in V : r(y_\ell^*) \cap r(z_i) \not= \emptyset\}} \alpha_i^2 \biggr)^\frac{1}{2} \\
    & \leq \frac{2}{m^2_{2j+1}}.
\end{split}
\end{equation*}

\n where the first inequality follows by applying Lemma \ref{lem:sw}
for each $i$ with ``$x^*$''$=y^*_i$, noting that its hypothesis is
satisfied by our assumption (1), and for the second inequality we
used the Cauchy-Schwartz inequality.

\end{proof}

{\footnotesize
\noindent
Department of Mathematics, University of South Carolina, Columbia, SC 29208, \\
\noindent giorgis@math.sc.edu, beanland@math.sc.edu
}

\end{document}